\documentclass[twoside]{article}
\usepackage{color}
\usepackage{latexsym}
\usepackage{amsmath}
\usepackage{amsfonts}
\usepackage{amssymb}
\usepackage{graphicx}
\usepackage{indentfirst}
\usepackage{makeidx}
\usepackage{url}
\usepackage{wrapfig}
\usepackage[colorlinks=true]{hyperref}
\usepackage{times,url}
\usepackage[normal]{subfigure}

\begin{document}

\def\reals{\mathbb R}
\def\wolp{\accentset{\circ}{W}_p^{1}}
\def\wlp{W_p^{1}}
\newcommand{\fin}{\hfill\mbox{$\quad{}_{\Box}$}}
\newcommand{\fineq}{\vspace{-.75cm$\fin$}\par\bigskip}
\newcommand{\fineqnum}{\vspace{-.4cm$\fin$}\par\bigskip}

\newtheorem{thm}{Theorem}[section]
\newtheorem{lem}{Lemma}[section]
\newtheorem{cor}{Corollary}[section]
\newtheorem{claim}{Proposition}[section]
\newtheorem{rem}{Remark}[section]
\newtheorem{defi}{Definition}[section]
\newtheorem{example}{Example}

\title{\Large \bfseries \sffamily On the exact multiplicity of stable ground states of non-Lipschitz
semilinear elliptic equations for some classes of starshaped sets}
\author{\bfseries\sffamily J.I. D\'{\i}az, J.~Hern\'{a}ndez and Y.Sh.~Ilyasov }
\date{}

\maketitle
\begin{abstract}
We prove the exact multiplicity of flat and compact support
stable solutions of an autonomous non-Lipschitz semilinear elliptic
equation of eigenvalue type according to the dimension N
and the two exponents, $0<\alpha<\beta<1$, of the involved nonlinearites. Suitable assumptions are made on the spatial domain $\Omega$ where the problem is
formulated in order to avoid a possible continuum of those solutions
and, on the contrary, to ensure the exact number of solutions
according to the nature of the domain $\Omega$. Our results also clarify some
previous works in the literature. The main techniques of proof are a
Pohozhaev's type identity and some fibering type arguments in the
variational approach.
\end{abstract}

\section{Introduction}

In this paper we study the existence of non-negative solutions of the
following problem
\begin{equation}  \label{PL}
\begin{cases}
-\Delta u+|u|^{\alpha -1}u=\lambda |u|^{\beta -1}u~~\mbox{in}~\Omega , \\
~~u=0~~\mbox{on}~\partial \Omega .%
\end{cases}
\tag*{$P(\alpha ,\beta ,\lambda )$}
\end{equation}%
Here $\Omega $ is a bounded domain in $\mathbb{R}^{N}$, $N\geq 3$ with a
smooth boundary $\partial \Omega $, which is strictly star-shaped with
respect to a point $x_{0}\in $ $\mathbb{R}^{N}$ (which will be identified as
the origin of coordinates if no confusion may arise), $\lambda $ is a real
parameter, $0<\alpha <\beta <1$. By a weak solution of $P(\alpha ,\beta
,\lambda ) $ we mean a critical point $u\in H_{0}^{1}:=H_{0}^{1}(\Omega ) $
of the energy functional
\begin{equation*}
E_{\lambda }(u)=\frac{1}{2}\int_{\Omega }|\nabla u|^{2}dx+\frac{1}{{\alpha +1%
}}\int_{\Omega }|u|^{{\alpha +1}}dx-\frac{\lambda }{\beta +1}\int_{\Omega
}|u|^{\beta +1}dx,
\end{equation*}%
where $H_{0}^{1}(\Omega )$ is the standard vanishing on the boundary Sobolev
space. We are interested in \textit{ground states } of $P(\alpha ,\beta
,\lambda )$: i.e., a weak solution $u_{\lambda }$ of $P(\alpha ,\beta
,\lambda )$ which satisfies the inequality
\begin{equation*}
E_{\lambda }(u_{\lambda })\leq E_{\lambda }(w_{\lambda })
\end{equation*}%
for any non-zero weak solution $w_{\lambda }$ of $P(\alpha ,\beta ,\lambda )$%
. Notice that in \cite{PucciSerrin} authors also use
the term \textquotedblleft ground state\textquotedblright\ with a different
meaning.

Since the diffusion-reaction balance $-\Delta u=f(\lambda ,u)$ involves the
non-linear reaction term
\begin{equation*}
f(\lambda ,u):=\lambda |u|^{\beta -1}u-|u|^{\alpha -1}u,
\end{equation*}%
and it is a non-Lipschitz function at zero (since $\alpha <1$ and $\beta <1$%
) important peculiar behavior of solutions of these problems arises. For
instance, that may lead to the violation of the Hopf maximum principle on
the boundary and the existence of compactly supported solutions as well as
the so called \ \textit{flat solutions} which correspond to weak solutions $%
u>0$ in $\Omega $ such that
\begin{equation}
\frac{\partial u}{\partial \nu }=0~~\mbox{on}~~\partial \Omega ,  \label{N}
\end{equation}%
where $\nu $ denotes the unit outward normal to $\partial \Omega $. When the
additional information (\ref{N}) holds but the weak solution may vanish in a
positively measured subset of $\Omega $, i.e. if $u\geq 0$ in $\Omega $, we
shall call it as a \textit{compact support solution} of $P(\alpha ,\beta
,\lambda )$ (sometimes also called as a \textit{free boundary solution, }%
since the boundary of its support is not a a priori known). Notice that in
that case the support of $u$ is strictly included in $\overline{\Omega }$.
If $u$ is a weak solution such that property (\ref{N}) is not satisfied we
shall call it as a \textit{usual weak solution} (since, at least for the
associated linear problem and for Lipschitz non-linear terms, the strong
maximum principle due to Hopf, implies that (\ref{N}) cannot be verified).

In what follows we shall use the following notation: any largest ball $%
B_{R(\Omega )}:=\{x\in \mathbb{R}^{N}:~|x|\leq R(\Omega )\}$ contained in $%
\Omega $ will be denoted as an \textit{inscribed ball} in $\Omega $. Our
exact multiplicity results will concern the case of some classes of
starshaped sets of $\mathbb{R}^{N}$ containing a finite number of different
\textit{inscribed balls} in $\Omega .$

For sufficiently large $\lambda $ the existence of a \textit{compactly
supported solution} of $P(\alpha ,\beta ,\lambda )$ follows from \cite%
{GazzolaSerin, Serrin-Zou} (see also for the case $N=1$, \cite{diaz,
Diaz-Hernan-Man}, \cite{CortElgFelmer-1, CortElgFelmer-2, Kaper1, Kaper2}.
Indeed, by \cite{Kaper1, Kaper2, GazzolaSerin, Serrin-Zou} the equation in $%
P(\alpha ,\beta ,1)$ considered in $\mathbb{R}^{N}$ has a unique (up to
translation in $\mathbb{R}^{N}$) compactly supported solution $u^{\ast }$, moreover $%
u^{\ast }$ is radially symmetric such that supp$(u^{\ast })$=$\overline{B}%
_{R^{\ast }}$ for some $R^{\ast }>0$. Hence since the support of $u_{\sigma
}^{\ast }(x):=u^{\ast }(x/\sigma ),~x\in B_{\sigma R^{\ast }}$ is contained
in $\Omega ,$ for sufficiently small $\sigma $, the function $w_{\lambda
}^{c}(x)=\sigma ^{-\frac{2}{1-\alpha }}\cdot u_{\sigma }^{\ast }(x)$ weakly
satisfies $P(\alpha ,\beta ,\lambda )$ in $\Omega $ with $\lambda =\sigma ^{-%
\frac{2(\beta -\alpha )}{1-\alpha }}$. However, it is not hard to show (see,
e.g. Corollary 5.2 below, that, in general (for all sufficiently large $%
\lambda $), weak solutions $w_{\lambda }$ are not ground states.

On the other side, finding \textit{flat} or \textit{compactly
supported ground states} is important in view of the study of
non-stationary problems (see \cite{DIH1, DIH, ilCrit} and
\cite{Rosenau}).

The existence of \textit{flat }and \textit{compact support ground states},
for certain $\lambda ^{\ast }$ of $P(\alpha ,\beta ,\lambda )$ has been
obtained in \cite{IlEg} (see also \cite{DIH}). In the present paper we
develop this result presenting here a sharper explanation of the main
arguments of its proof. Furthermore, we shall offer here some more precise
results on the behaviour of ground states depending on $\lambda $.

It is well known that the non-Lipschitz nonlinearities may entail
the existence of a continuum of nonnegative compact supported
solutions of elliptic boundary value problems. However the answer
for the same question stated about ground states or \textit{usual}
solutions becomes unclear. Notice that this question is important in
the investigation of stability solutions for non-stationary problems
(see \cite{DIH1, DIH, ilCrit}). We recall that, as a matter of fact,
flat solutions of  $P(\alpha ,\beta ,\lambda ^{\ast })$ only may
arise if $\Omega $ is the ball $B_{R^{\ast }}$ mentioned before. For
the rest of domains, and values of $\lambda \geq \lambda ^{\ast },$
any weak solution which is not a \textquotedblleft
usual\textquotedblright\ solution should have compact support.

%This rise a question {\it as to whether the similar phenomena may be
%occurred in case of the higher dimensions $n>1$}. More precisely
%whether the Hopf boundary maximum principle  holds for
%(\ref{L}) when $n>1$ and  the nonlinearity
%$f(\lambda,u)$ is non-Lipschitz. To find an answer to this
%question is a main goal in the present work.
%
%A central role plays the following
%\begin{equation}
%d^*=\frac{(N-p)}{Np}\alpha\beta +2-\beta-\alpha.
%\end{equation}

%Let us state our main result.
%We consider a weak solution $u \in W:=W(\Omega)$, where $W(\Omega)$ denotes the closure $C^\infty_0(\Omega)$ in standard Sobolev space $H^1(\Omega)$ with the norm $||\cdot||_1$. We say that a weak solution $u \in W$ of (\ref{L}) is a non-regular, if $u \in C^1(\overline{\Omega})$ and $\frac{\partial u}{\partial \nu} =0$ on $\partial \Omega$.

Let us state our main results. For given $u\in H_{0}^{1}(\Omega )$,
the \textit{fibrering mappings} are defined by $\phi
_{u}(t)=E_{\lambda }(tu)$ so that from the variational formulation
of $P(\alpha ,\beta ,\lambda )$ we know that $\phi _{u}^{\prime
}(t)|_{t=1}=0$ for solutions, where we use the notation
\begin{equation*}
\phi _{u}^{\prime }(t)=\frac{\partial }{\partial t}E_{\lambda }(tu).
\end{equation*}%
%
%
%
%
%
%
%
%
%
%
%
%
%
%
%
%
%
%\begin{figure}[th]
%\center{\includegraphics[width=0.55\linewidth]{fig2.pdf}}
%\caption{$r_{\tiny min}$ and $r_{\tiny max}$}
%\end{figure}
If we also define $\phi _{u}^{\prime \prime }(t)=\frac{\partial ^{2}}{%
\partial t^{2}}E_{\lambda }(tu)$, then, in case $\beta <1$ the equation $%
\phi _{u}^{\prime }(t)=0$ may have at most two nonzero roots $t_{\min }(u)>0$
and {$t_{\max }(u)>0$} such that $\phi _{u}^{\prime \prime }(t_{\max
}(u))\leq 0$, $\phi _{u}^{\prime \prime }(t_{\min }(u))\geq 0$ and $%
0<t_{\max }(u)\leq t_{\min }(u)$. This implies that any weak solution of $%
P(\alpha ,\beta ,\lambda )$ (any critical point of $E_{\lambda }(u)$)
corresponds to one of the cases $t_{\min }(u)=1$ or $t_{\max }(u)=1$.
However, it was discovered in \cite{IlEg} (see also \cite{DIH, ilDr, ilCrit}%
) that in case when we study flat \ or compactly supported solutions this
correspondence essentially depends on the relation between $\alpha $, $\beta
$ and $N$. Thus following this idea (from \cite{DIH, ilDr, ilCrit, IlEg}, in
the case $N\geq 3$, we consider the following subset of exponents
\begin{equation*}
\mathcal{E}_{s}(N):=\{(\alpha ,\beta ):~~2(1+\alpha )(1+\beta )-N(1-\alpha
)(1-\beta )<0,~0<\alpha <\beta <1\}.
\end{equation*}%
The main property of $\mathcal{E}_{s}(N)$ is that for star-shaped domains $%
\Omega $ in $\mathbb{R}^{N}$, $N\geq 3,$ if $(\alpha ,\beta )\in \mathcal{E}%
_{s}(N)$, any ground state solution $u$ of $P(\alpha ,\beta ,\lambda )$
satisfies $\phi _{u}^{\prime \prime }(t)|_{t=1}>0$ (see Lemma \ref{pro}
below and \cite{DIH, IlEg}).

\begin{rem}
In the cases $N=1,2$, one has $\mathcal{E}_{s}(N)=\emptyset $. Furthermore,
this implies (see \cite{DIH}) that if $N=1,2$ and $0<\alpha <\beta <1$, then
any flat or compact support weak solution $u$ of $P(\alpha ,\beta ,\lambda )$
satisfies $\phi _{u}^{\prime \prime }(t)|_{t=1}<0$.
\end{rem}

%We need also in the following  extremal values  introduced in \cite{DIH, IlEg} ( which are in fact the nonlinear generalized Rayleigh extremal values (see \cite{ilyaReil} and below for details))
%\begin{equation}\label{PNR}
%\Lambda_0=c_0^{\alpha,\beta}\cdot\Lambda, ~~\Lambda_1=c_1^{\alpha,\beta}\cdot\Lambda,
%\end{equation}
%where
%\begin{equation}
%\Lambda=\inf_{W\setminus 0}\frac{(\int |u|^{{\alpha+1}} dx)^\frac{1-\beta}{1-\alpha}(\int|\nabla u|^{2} dx)^\frac{\beta-\alpha}{1-\alpha}}{\int |u|^{\beta+1} dx}
%\end{equation}
%and $c_1^{\alpha,\beta},c_0^{\alpha,\beta}$ are constants which are expressed by explicit formulas \eqref{cE}, \eqref{cPD} and satisfy $c_1^{\alpha,\beta}<c_0^{\alpha,\beta}$.
%

In what follows we shall use the notations
\begin{equation*}
E_{\lambda }^{\prime }(u)=\phi _{u}^{\prime }(t)|_{t=1}=\frac{\partial }{%
\partial t}E_{\lambda }(tu)|_{t=1},~~E_{\lambda }^{\prime \prime }(u)=\phi
_{u}^{\prime \prime }(t)|_{t=1}=\frac{\partial ^{2}}{\partial t^{2}}%
E_{\lambda }(tu)|_{t=1},~~~u\in H_{0}^{1}(\Omega ).
\end{equation*}%
Our first result is the following %\begin{thm}\label{Th1}

\begin{thm}
\label{Th1} Let $N\geq 3$ and let $\Omega $ be a bounded strictly
star-shaped domain in $\mathbb{R}^{N}$ with $C^{2}$-manifold boundary $%
\partial \Omega $. Assume that $(\alpha ,\beta )\in \mathcal{E}_{s}(N)$.
Then there exists $\lambda ^{\ast }>0$ such that for any $\lambda \geq
\lambda ^{\ast }$ problem $P(\alpha ,\beta ,\lambda )$ possess a ground
state $u_{\lambda }$. Moreover $E_{\lambda }^{\prime \prime }(u_{\lambda
})>0 $, $u_{\lambda }\in C^{1,\gamma }(\overline{\Omega })$ for some $\gamma
\in (0,1)$ and $u_{\lambda }\geq 0$ in $\Omega $. For any $\lambda <\lambda
^{\ast }$, problem $P(\alpha ,\beta ,\lambda )$ has no weak solution.
\end{thm}

\bigskip

Our second main result deals with the (non-)existence of flat or compactly
supported ground states.

\begin{thm}
\label{Th2} There is a non-negative ground state $u_{\lambda ^{\ast }}$
which is flat or has compact support. Moreover, $u_{\lambda ^{\ast }}$ is
radially symmetric about some point of $\Omega $, and supp$(u_{\lambda
^{\ast }})$=$\overline{B}_{R(\Omega )}$ is an inscribed ball in $\Omega $.
For all $\lambda >\lambda ^{\ast }$, any ground state $u_{\lambda }$ of $%
P(\alpha ,\beta ,\lambda )$ is a \textquotedblleft usual\textquotedblright\
solution.
\end{thm}

\bigskip

Our last result deals with the multiplicity of solutions. Our main goal is
to extend the results of \cite{diaz} and \cite{Diaz-Hernan-Man} concerning
the one-dimensional case. We also recall that the existence of what we call
now \textquotedblleft usual\textquotedblright\ solutions was proved in some
previous papers in the literature. Existence of a smooth branch of such
positive solutions was proved for $\lambda >\lambda ^{\ast }$ in \cite{HMV}
by using a change of variables and then a continuation argument. The
existence of at least two non-negative solutions in such a case was shown in
\cite{Montenegro} by using variational arguments and this result was
improved in \cite{Annello} showing that one of the solutions is actually
positive, again by variational arguments. Many of these results are valid
even in the singular case $-1<\alpha <\beta <1.$

In order to present our exact multiplicity results we introduce the \textit{%
geometrical reflection} across a given hyperplane $H$ by the usual isometry $%
R_{H}:\mathbb{R}^{N}\rightarrow $ $\mathbb{R}^{N}$. Remember that any point
of $H$ is a fixed point of $R_{H}$. Now we shall introduce some classes of
starshaped sets $\Omega $ for which we can obtain the exact multiplicity of
flat stable ground solutions of problem $P(\alpha ,\beta ,\lambda ^{\ast }).$
We say that $\Omega $ is of \textit{Strictly} \textit{Starshaped Class }$%
\mathit{m}$, if it is a strictly starshaped domain and contains
exactly $m$ inscribed balls of the same radius $R(\Omega )$ such
that each of them can be obtained from any other by $k\in
\{1,...,m\}$ reflections of $\Omega $ across some hyperplanes
$H_{i}$, $i=1,...,k$.

\bigskip

\begin{thm}
\label{ThmCor1} Assume $N\geq 3$, $(\alpha ,\beta )\in \mathcal{E}_{s}(N)$.
Let $\Omega $ be a domain of Strictly Starshaped Class $m>1$ with a $C^{2}$%
-manifold boundary $\partial \Omega $. Then there exist exactly $m$
stable nonnegative flat or compact supported ground states
$u_{\lambda ^{\ast }}^{1}$, $u_{\lambda ^{\ast }}^{2}$,...,
$u_{\lambda ^{\ast }}^{m}$ of problem $P(\alpha ,\beta ,\lambda
^{\ast })$ and $m$ sets of
\textquotedblleft \textit{usual\textquotedblright } ground states $%
(u_{\lambda _{n}}^{1})_{n=1}^{\infty }$, $(u_{\lambda
_{n}}^{2})_{n=1}^{\infty }$,..., $(u_{\lambda _{n}}^{m})_{n=1}^{\infty }$ of
$P(\alpha ,\beta ,\lambda _{n}),$ , with $\lim_{n\rightarrow \infty }\lambda
_{n}=\lambda ^{\ast }$, $\lambda _{n}>\lambda ^{\ast }$, $n=1,2,...$ and
such that $u_{\lambda _{n}}^{i}\rightarrow u_{\lambda ^{\ast }}^{i}$,
strongly in $H_{0}^{1}$ as $n\rightarrow \infty ,$ for any $i=1,...,m$.
\end{thm}

%\textcolor[rgb]{1,0,0}{I arose in the statement of theorem the word {\it flat }, i.e. it  is written now  ...\textcolor[rgb]{0,0,1}{ compact support stable nonnegative ground
%solutions of problem ...}. This has been done since the solution  can not be flat in the case $m>1$.}

Let us show how can be obtained some domains of Strictly Starshaped class $m$%
. We start by considering an initial bounded Lipschitz set $\Omega _{1}$ of $%
\mathbb{R}^{N}$ such that:
\begin{equation}
\Omega _{1}\text{ contains exactly one inscribed ball of radius }R(\Omega_{1} )
\text{.}  \label{First cond omega}
\end{equation}%
We also introduce the following notation: given a general open set $G$ of $%
\mathbb{R}^{N}$ we define $S[G]$ as the set of points $y\in G$ such that $G$
is strictly starshaped with respect to $y$. Then, the second condition we
shall require to $\Omega _{1}$ is
\begin{equation}
S[\Omega _{1}]~~\text{ is not empty.}  \label{starshap}
\end{equation}%
Then $\Omega $ belongs to the Strict Starshaped class $1$ if there
exists $\Omega
_{1}$ satisfying (\ref{First cond omega}) and (\ref{starshap}) such that $%
\Omega =\Omega _{1}$. Now, let us show how we can obtain a domain of
Strictly Starshaped class $2$.

Let $\Omega _{1}$ be a domain of Strictly Starshaped class $1$ and assume,
additionally, that the set $S[\Omega _{1}]$ contains some other point
different than $x_{1}$, $\{x_{1}\}\varsubsetneq S[\Omega _{1}]$, i.e.%
\begin{equation*}
\text{there exists }y_{1}\in S[\Omega _{1}]\text{ such that }y_{1}\neq x_{1}.
%\label{not centered}
\end{equation*}%
Let now $\Omega _{2}:=R_{H(y_{1})}(\Omega _{1})$ be the reflected set of $%
\Omega _{1}$ across some hyperplane $H(y_{1})$ containing the point $y_{1}$
such that
\begin{equation*}
\Omega _{1}\cup \Omega _{2}\text{ contains exactly one inscribed ball of
radius }R(\Omega)\text{ of center }x_{2}\neq x_{1}%
\text{.}
\end{equation*}%
We now consider
\begin{equation*}
\Omega =\Omega _{1}\cup \Omega _{2}.
\end{equation*}%
Notice that, obviously, $\Omega $ is Strictly Starshaped class $1$ with
respect to $y_{1}$ (since $y_{1}\in S[\Omega _{1}]$ and any ray starting
from $y_{1}$ is reflected to a ray linking $y_{1}$ with any other point of $%
\Omega _{2}$). Moreover, such a domain $\Omega $ verifies
\begin{align*}
\Omega &\text{ contains exactly two inscribed balls of radius }{R(\Omega )},\\
&\text{ with center at two different points }x_{i}\in \Omega ,\text{ }i=1,2.
\end{align*}%
Thus $\Omega $ is a set of \textit{Strictly Starshaped class 2}. Evidently
we can repeat this construction with a domain of \textit{Strictly Starshaped
class 2} and obtained a domain $\Omega $ of \textit{Strictly Starshaped
class 3}, etc.

\begin{figure}[!ht]
\begin{center}
\center{\includegraphics[scale=0.4]{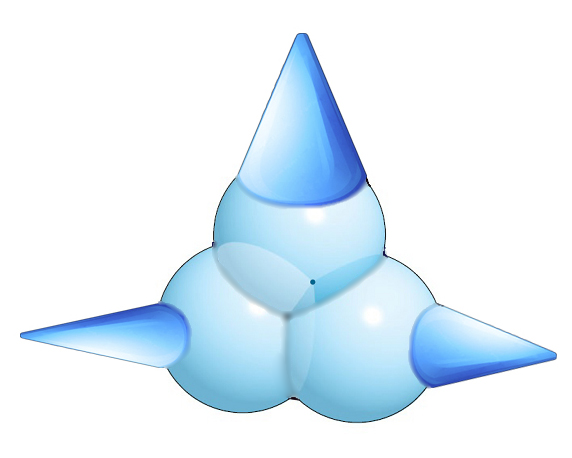}}\\[0pt]
\end{center}
\caption{Domain generating exactly three ground states}
\end{figure}

We believe that we can iterate this process in a similar way until
some number $m:=m(N)\geq 3$, which maybe depends on the dimension
$N$. However we don't know how to prove this. Moreover we rise the
following conjecture: \textit{For a given dimension }$N$, \textit{\
there exists a number $m(N)$ such that for any $k=1,2,...,m(N)$
there exists a domain of Strictly Starshaped class $k$ whereas there
is no domain in $\mathbb{R}^{N}$ of Strictly Starshaped class $k$
with $k>m(N)$.}

\begin{figure}[!ht]
\begin{center}
\center{\includegraphics[scale=0.4]{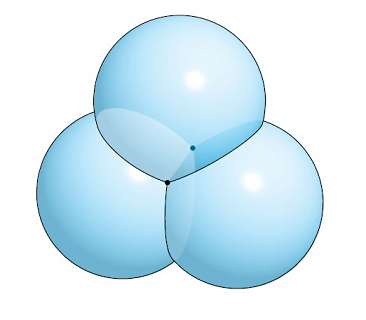}}\\[0pt]
\end{center}
\caption{Union of the supports of the three radially symmetric ground states
corresponding to the domain given by Figure 1.}
\end{figure}

\begin{rem}
We emphasize that by Theorems \ref{Th1}, \ref{Th2},\ref{ThmCor1} we
obtain the complete bifurcation diagram for the ground states of
$P(\alpha ,\beta ,\lambda )$ for domains of Starshaped Class $m$.
Indeed, the flat ground state $u_{\lambda ^{\ast }}$ corresponds to
a fold bifurcation point (or turning point) from which it start
$m+1$ different branches of weak solutions: on one hand, the branch
of \textquotedblleft usual\textquotedblright\ ground states
$u_{\lambda }$, forming a branch of stable equilibria, and, on the
other hand, $m$ branches formed by unstable compactly supported weak
solutions, of the form $w_{\lambda }^{c}(x:x_{0,j})=\sigma
^{-\frac{2}{1-\alpha }}\cdot u_{\lambda ^{\ast
}}((x-x_{0,j})/\sigma )$ with $\lambda =\sigma ^{-\frac{2(\beta -\alpha )}{%
1-\alpha }}$ (see Figure 1) and $m$ different points $x_{0,j}$,
$j=1,...,m$. Furthermore, we know a global information: the energy
of $u_{\lambda ^{\ast }}$ is the
maximum among all the possible energies associated to any weak solution of $%
P(\alpha ,\beta ,\lambda )$.
\end{rem}

\begin{figure}[th]
\begin{center}
\center{\includegraphics[scale=0.4]{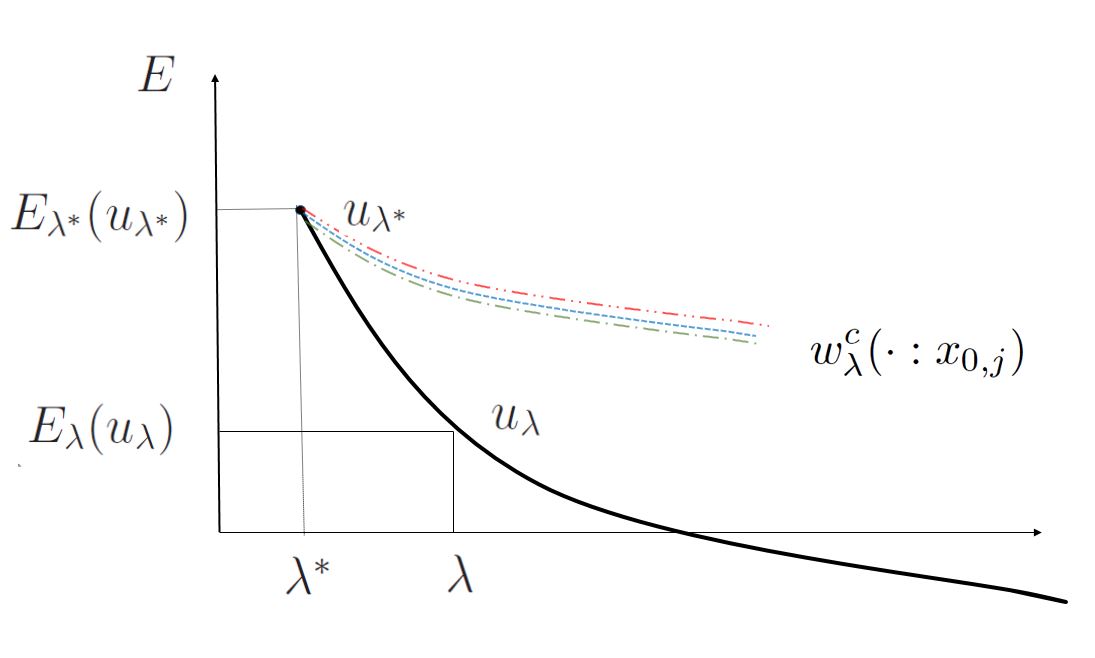}}\\[0pt]
\end{center}
\caption{Bifurcation diagram for the energy levels of ground states
and compact support solutions. }
\end{figure}
In the last part of the paper we consider the associate parabolic problem

\begin{equation}
PP(\alpha ,\beta ,\lambda ,v_{0})\quad \left\{
\begin{array}{ll}
v_{t}-\Delta v+|v|^{\alpha -1}v=\lambda |v|^{\beta -1}v & \text{in }%
(0,+\infty )\times \Omega \\
v=0 & \text{on }(0,+\infty )\times \partial \Omega \\
v(0,x)=v_{0}(x) & \text{on }\Omega .%
\end{array}%
\right.  \label{p1}
\end{equation}%
For the basic theory for this problem, under the structural assumption $%
0<\alpha <\beta <1$ we send the reader to \cite{DIH} and its references. We
apply here some local energy methods, for the two cases $\lambda >\lambda
^{\ast }$ and $\lambda =\lambda ^{\ast }$, to give some information on the
evolution and formation, respectively, of the free boundary given by the
boundary of the support of the solution $v(t,.)$ when $t$ increases. This
provides a complemmentary information since by Theorem 1.1 (and the
asymptotic behaviour results for $PP(\alpha ,\beta ,\lambda ,v_{0})$) we
know that, as $t\rightarrow +\infty $, the support of $v(t,.)$ must converge
to a ball of $\mathbb{R}^{N}$, in the case $\lambda =\lambda ^{\ast }$, or
to the whole domain $\overline{\Omega }$, if $\lambda >\lambda ^{\ast }$,
(the supports of one of the corresponding stationary solutions).

%The proof of the theorem relies on the variational arguments.
%Furthermore, basic ingredients in the proof consist in using  Pohozaev's
%identity \cite{poh} corresponding to (\ref{L}) and in applying the spectral analysis with respect to the fibering procedure introduced in \cite{ilyas}.

%\begin{center}
%\bigskip In the following picture we show the union of the support of
%solutions associated to a domain $\Omega $ of Class $3$ \textit{starting
%from a the set }$\Omega _{1}$\textit{\ }$\FRAME{itbpF}{4.4356in}{4.1451in}{%
%0in}{}{\Qlb{union of the support of solutions associated to a domain O of
%Class 3 }}{three spheres.eps}{\special{language "Scientific Word";type
%"GRAPHIC";maintain-aspect-ratio TRUE;display "USEDEF";valid_file "F";width
%4.4356in;height 4.1451in;depth 0in;original-width 13.1158in;original-height
%12.2509in;cropleft "0";croptop "1";cropright "1";cropbottom "0";filename
%'three spheres.eps';file-properties "XNPEU";}}$
%\end{center}

\section{Preliminaries}

In this section we give some preliminary results. In what follows $%
H_{0}^{1}:=H_{0}^{1}(\Omega )$ denotes the standard vanishing on the
boundary Sobolev space. We can assume that its norm is given by
\begin{equation*}
||u||_{1}=\left( \int_{\Omega }|\nabla u|^{2}\,dx\right) ^{1/2}.
\end{equation*}%
Denote
\begin{equation*}
P_{\lambda }(u):=\frac{1}{2^{\ast }}\int_{\Omega }|\nabla u|^{2}\,\mathrm{d}%
x+\frac{1}{{\alpha +1}}\int_{\Omega }|u|^{{\alpha
+1}}\,\mathrm{d}x-\lambda \frac{1}{{\beta +1}}\int_{\Omega
}|u|^{{\beta +1}}\,\mathrm{d}x,
\end{equation*}%
where%
\begin{equation*}
2^{\ast }=\frac{2N}{N-2}~~\text{ for }N\geq 3.
\end{equation*}%
We will use the notation $P_{\lambda }^{\prime }(tu)=dP_{\lambda }(tu)/dt$, $%
t>0$, $u\in H_{0}^{1}$. From now on we suppose that the boundary $\partial
\Omega $ is a $C^{2}$-manifold. As usual, we denote by $\mathrm{d}\sigma $
the surface measure on $\partial \Omega $. We need the Pohozhaev's identity
for a weak solution of $P(\alpha ,\beta ,\lambda ).$

\begin{lem}
\label{lem1} Assume that $\partial \Omega $ is a $C^{2}$-manifold, $N\geq 3$%
. Let $u\in C^{1}(\overline{\Omega })$ be a weak solution of $P(\alpha
,\beta ,\lambda )$. Then there holds the Pohozaev identity
\begin{equation*}
P_{\lambda }(u)=-\frac{1}{2N}\int_{\partial \Omega }\left\vert \frac{%
\partial u}{\partial \nu }\right\vert ^{2}\,(x\cdot \nu (x))\mathrm{d}\sigma
(x).
\end{equation*}
\end{lem}

For the proof see \cite{DIH, Takac_Ilyasov}and \cite{poh},
\cite{Temam}. See also some related results in \cite{poh} and
\cite{Temam}.

\noindent Notice that
\begin{equation}
E_{\lambda }(u)=P_{\lambda }(u)+\frac{1}{N}\int_{\Omega }|\nabla
u|^{2}dx,~~~\forall u\in H_{0}^{1}(\Omega ).  \label{PandE}
\end{equation}

Assume $\Omega $ is strictly star-shaped with
respect to a point $x_{0}\in $ $\mathbb{R}^{N}$ (which will be identified as
the origin of coordinates of $\mathbb{R}^{N}$).  Observe that if $\Omega $
is a star-shaped (strictly star-shaped) domain with respect to the origin of
$\mathbb{R}^{N}$, then $x\cdot \nu \geq 0$ ($x\cdot \nu >0$) for all $x\in
\partial \Omega $. This and Lemma \ref{lem1} imply

\begin{cor}
\label{cor1} Let $\Omega $ be a bounded star-shaped domain in $\mathbb{R}%
^{N} $ with a $C^{2}$-manifold boundary $\partial \Omega $. Then any weak
solution $u\in C^{1}(\overline{\Omega })$ of $P(\alpha ,\beta ,\lambda )$
satisfies $P_{\lambda }(u)\leq 0$. Moreover, if $u$ is a flat solution or it
has a compact support then $P_{\lambda }(u)=0$. Furthermore, in the case $%
\Omega $ is strictly star-shaped, the converse is also true: if $P_{\lambda
}(u)=0$ and $u\in C^{1}(\overline{\Omega })$ is a weak solution of $P(\alpha
,\beta ,\lambda )$, then $u$ is flat or it has a compact support.
\end{cor}

The proof of the following result can be found in \cite{DIH, IlEg}.

\begin{lem}
\label{pro} Assume $N\geq 3$ and $(\alpha ,\beta )\in \mathcal{E}_{s}(N)$.

$(i)$ Let $u\in C^{1}(\overline{\Omega })$ be a flat or compact support weak
solution of $P(\alpha ,\beta ,\lambda )$. Then $E_{\lambda }(u)>0$ and $%
E_{\lambda }^{\prime \prime }(u)>0$.

$(ii)$ If $E_{\lambda }^{\prime }(u)=0$, $P_{\lambda }(u)\leq 0$ for some $%
u\in H_{0}^{1}(\Omega )\setminus 0$, then
\begin{equation*}
E_{\lambda }^{\prime \prime }(u)>0.
\end{equation*}
\end{lem}

\begin{rem}
When $0<\beta <\alpha <1,$ a case which is not considered in this paper, we
have $E_{\lambda }^{\prime \prime }(u)>0$ and $P_{\lambda }(u)<0$ for any
weak solution $u\in H_{0}^{1}\setminus 0$ of $P(\alpha ,\beta ,\lambda )$.
In particular, in this case, any solution of $P(\alpha ,\beta ,\lambda )$ is
a \textquotedblleft usual\textquotedblright\ solution. The uniqueness of
solution was shown in \cite{HMV}.
\end{rem}

In what follows we need also

\begin{claim}
\label{pradd} If $E_{\lambda}^{\prime }(tu)=0$ for $u\neq 0 $, then $%
P^{\prime }_\lambda(tu)<0$.
\end{claim}
\noindent {\em Proof}\quad Observe that,
\begin{equation*}
P_{\lambda }^{\prime }(tu)=\frac{N-2}{N}t\int_{\Omega }|\nabla u|^{2}\,%
\mathrm{d}x-\lambda t^{\beta }\int_{\Omega }|u|^{\beta +1}\,\mathrm{d}%
x+t^{\alpha }\int_{\Omega }|u|^{\alpha +1}\,\mathrm{d}x=E_{\lambda
}^{\prime }(tu)-\frac{2t}{N}\int_{\Omega }|\nabla
u|^{2}\,\mathrm{d}x.
\end{equation*}%
Thus $E_{\lambda }^{\prime }(tu)=0$ entails $P^{\prime}_{\lambda
}(tu)=-(2t/N)\int |\nabla u|^{2}\,\mathrm{d}x<0\fin$.

\section{Auxiliary extremal values}

In this section we introduce some extremal values which will play an
important role in the following. Some of these values, and the corresponding
variational functionals, have been already introduced in \cite{DIH, IlEg}.
However, for our aims we shall introduce them using another approach which
is more natural and easy.

Our approach will be based on using a \textit{nonlinear generalized Rayleigh
quotient} (see \cite{ilyaReil}). In fact, we can associate to problem $%
P(\alpha ,\beta ,\lambda )$ several nonlinear generalized Rayleigh quotients
which may give useful information on the nature of the problem. In this
paper we will deal with two of them.

First, let us consider the following Rayleigh's quotient \cite{ilyaReil}
\begin{equation}
R^{0}(u)=\frac{\frac{1}{2}\int_{\Omega }|\nabla u|^{2}dx+\frac{1}{{\alpha +1}%
}\int_{\Omega }|u|^{{\alpha +1}}dx}{\frac{1}{{\beta +1}}\int_{\Omega }|u|^{{%
\beta +1}}dx},~~u\neq 0.  \label{lamb1}
\end{equation}%
Following \cite{ilyaReil}, we consider
\begin{equation}
r_{u}^{0}(t):=R^{0}(tu)=\frac{\frac{t^{1-\beta }}{2}\int_{\Omega }|\nabla
u|^{2}dx+\frac{t^{\alpha -\beta }}{{\alpha +1}}\int_{\Omega }|u|^{{\alpha +1}%
}dx}{\frac{1}{{\beta +1}}\int_{\Omega }|u|^{{\beta +1}}dx},~~t>0,~~u\neq 0.
\label{RaylCC}
\end{equation}%
Notice that for any $u\neq 0$, \ and $\lambda \in \mathbb{R}$,
\begin{equation}
\text{if}~R^{0}(u)\equiv r_{u}^{0}(t)|_{t=1}=\lambda ,~~\text{then}%
~~E_{\lambda }(u)=0.  \label{LzeroProp}
\end{equation}%
%
%
%
%
%
%
%
%
%
%
%
%
%
%
%Compute
%~\text{and}~\frac{\partial}{\partial t} r_u^0(t)|_{t=1}=\frac{\beta+1}{\int |u|^{{\beta+1}} dx}E_\lambda^\prime(u)
%\begin{equation}
%\frac{\partial}{\partial t} r_u^0(t)=\frac{(1-\beta)\frac{t^{-\beta}}{2}\int|\nabla u|^{2} dx+(\alpha-\beta)\frac{t^{\alpha-\beta-1}}{{\alpha+1}}\int |u|^{{\alpha+1}} dx}{\frac{1}{\beta}\int |u|^{\beta} dx}.
%\end{equation}
It is easy to see that $\partial r_{u}^{0}(t)/\partial t=0$ if and only if
\begin{equation*}
(1-\beta )\frac{t^{-\beta }}{2}\int_{\Omega }|\nabla u|^{2}dx+(\alpha -\beta
)\frac{t^{\alpha -\beta -1}}{{\alpha +1}}\int_{\Omega }|u|^{{\alpha +1}}dx=0,
\end{equation*}%
and that the only solution to this equation is
\begin{equation}
t_{0}(u)=\left( \frac{2(\beta -\alpha )}{(\alpha +1)(1-\beta )}\frac{%
\int_{\Omega }|u|^{{\alpha +1}}dx}{\int_{\Omega }|\nabla u|^{2}dx}\right) ^{%
\frac{1}{1-\alpha }}.  \label{P11}
\end{equation}%
Let us emphasize that $t_{0}(u)$ is a value where the function $r_{u}^{0}(t)$
attains its global minimum. Substituting $t_{0}(u)$ into $r_{u}^{0}(t)$ we
obtain the following \textit{nonlinear generalized Rayleigh quotient}:
\begin{equation}
\lambda _{0}(u)=r_{u}^{0}(t_{0}(u))\equiv
R^{0}(tu)|_{t=t_{0}(u)}=c_{0}^{\alpha ,\beta }\lambda (u),
\label{lambdazero}
\end{equation}%
where
\begin{equation}
c_{0}^{\alpha ,\beta }=\frac{(1-\alpha )(\beta +1)}{(1-\beta )(1+\alpha )}%
\left( \frac{(1-\beta )(\alpha +1)}{2(\beta -\alpha )}\right) ^{\frac{\beta
-\alpha }{1-\alpha }},  \label{cE}
\end{equation}%
and
\begin{equation*}
\lambda (u)=\frac{(\int_{\Omega }|u|^{{\alpha +1}}dx)^{\frac{1-\beta }{%
1-\alpha }}(\int_{\Omega }|\nabla u|^{2}dx)^{\frac{\beta -\alpha }{1-\alpha }%
}}{\int_{\Omega }|u|^{{\beta +1}}dx}.
\end{equation*}%
See Figure 3.
%\FRAME{ftbpF}{7.1338in}{4.5558in}{0pt}{}{\Qlb{Curve lamnda}}{%
%figure1.eps}{\special{language "Scientific Word";type
%"GRAPHIC";maintain-aspect-ratio TRUE;display "USEDEF";valid_file "F";width
%7.1338in;height 4.5558in;depth 0pt;original-width 7.0681in;original-height
%4.5031in;cropleft "0";croptop "1";cropright "1";cropbottom "0";filename
%'Figure1.eps';file-properties "XNPEU";}}

\begin{figure}[!ht]
\begin{minipage}[t]{0.48\linewidth}
\center{\includegraphics[scale=0.4]{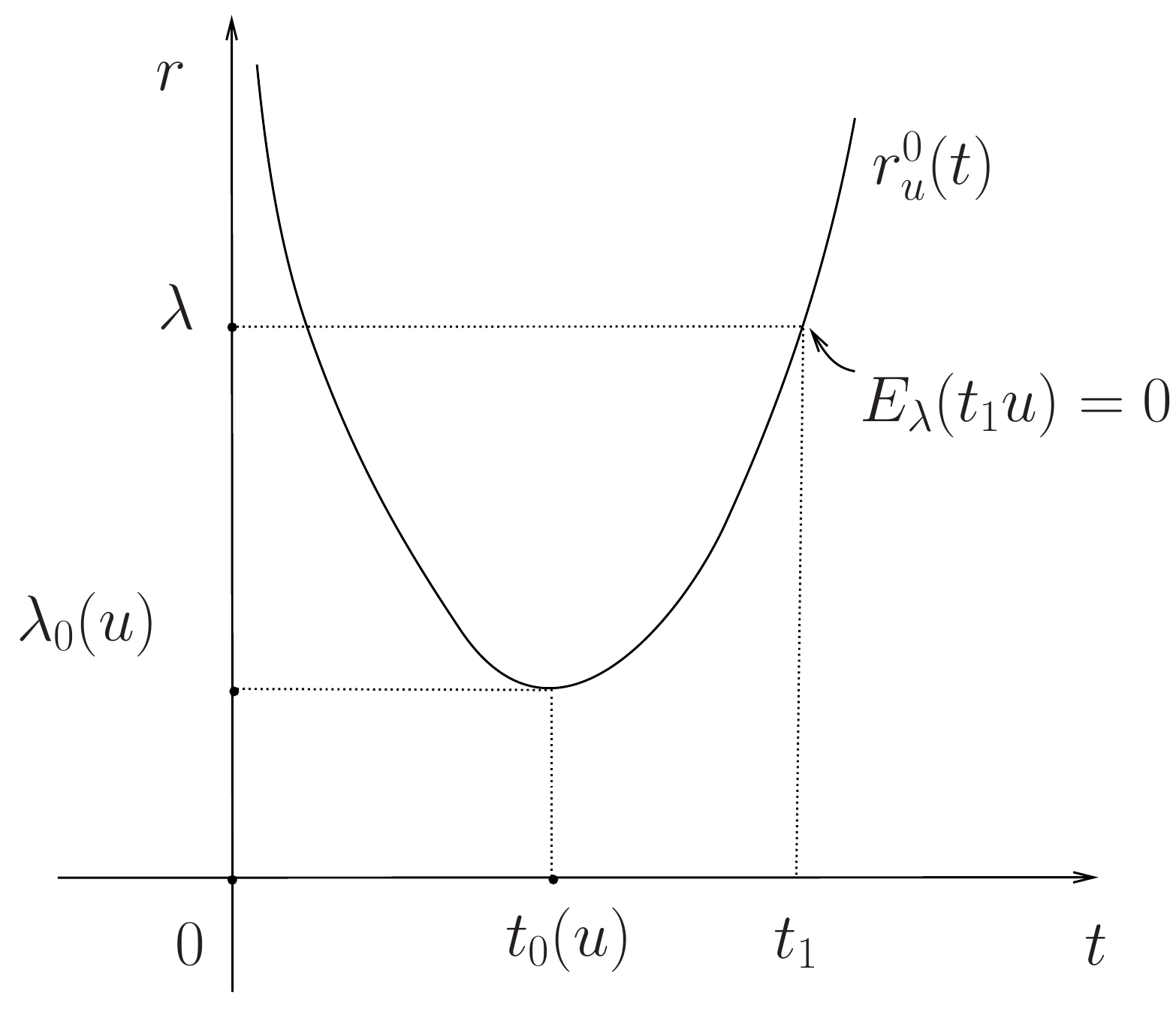}}\\
\caption{}
\end{minipage}
\hfill
\begin{minipage}[t]{0.48\linewidth}
\center{\includegraphics[scale=0.4]{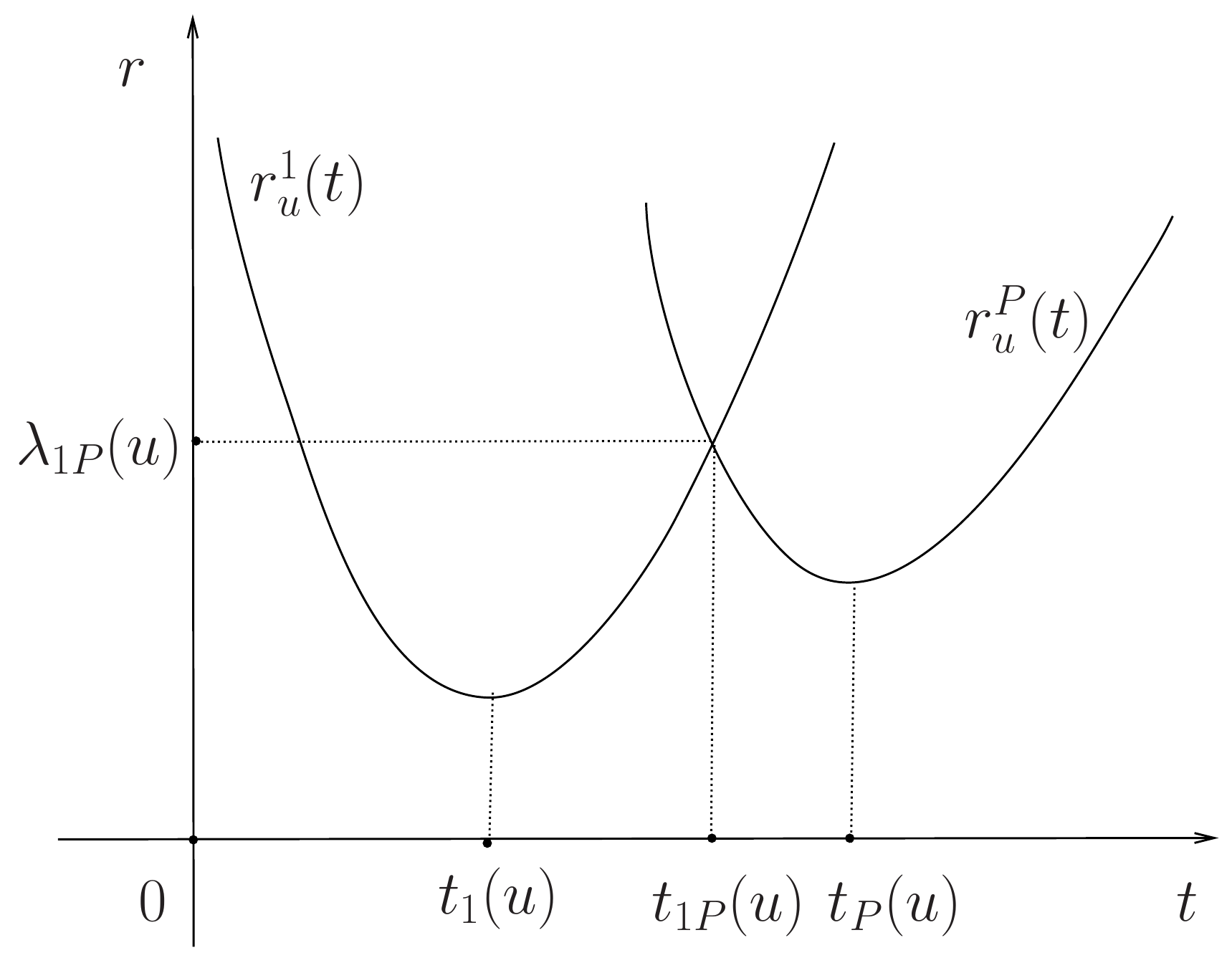}} \\
\caption{}
\end{minipage}
\end{figure}

It is not hard to prove (see, e.g., page 400 of \cite{Zeidler}) that

\begin{claim}
\label{propDiff} The map $\lambda (\cdot ):H_{0}^{1}(\Omega )\setminus
0\rightarrow \mathbb{R}$ is a $C^{1}$-functional.
\end{claim}

Consider the following extremal value of $\lambda _{0}(u)$
\begin{equation}
\Lambda _{0}=\inf_{u\in H_{0}^{1}(\Omega )\setminus 0}\lambda _{0}(u).
\label{P3}
\end{equation}%
Using Sobolev's and H\"{o}lder's inequalities (see, e.g.,
\cite{IlEg}) it can be shown that
\begin{equation}
0<\Lambda _{0}<+\infty .
\end{equation}%
By the above construction and using \eqref{LzeroProp} it is not hard to
prove the following

\begin{claim}
\label{propL0}

\begin{description}
\item[(i)] If $\lambda <\Lambda _{0}$, then $E_{\lambda }(u)>0$ for any $u
\neq 0$,

\item[(ii)] For any $\lambda >\Lambda _{0}$ there is $u\in H_{0}^{1}(\Omega
)\setminus 0$ such that $E_{\lambda }(u)<0$, $E_{\lambda }^{\prime }(u)=0$.
\end{description}
\end{claim}

%See Fig. 3
%\begin{figure}[!ht]
%\begin{minipage}[t]{0.48\linewidth}
%\center{\includegraphics[scale=0.4]{Fig3}}\\
%\caption{}
%\end{minipage}
%\end{figure}

In what follows we shall use the following result:

\begin{claim}
\label{CRRay} Let $u$ be a critical point of $\lambda _{0}(u)$ at some
critical value $\bar{\lambda}$, i.e. $D_{u}\lambda _{0}(u)=0 $, $\bar{\lambda%
}=\lambda _{0}(u)$. Then $D_{u}E_{\bar{\lambda}}(u)=0$ and $E_{\bar{\lambda}%
}(u)=0$.
\end{claim}

\noindent {\em Proof}\quad Observe that
\begin{equation*}
D_{u}\lambda _{0}(u)(\phi )=D_{u}r_{u}^{0}(t_{0}(u))(\phi )+\frac{\partial }{%
\partial t}r_{u}^{0}(t_{0}(u))(D_{u}t_{0}(u)(\phi ))=0,~~~~\forall \phi \in
C_{0}^{\infty }(\Omega ).
\end{equation*}%
Hence, since $\partial r_{u}^{0}(t)/\partial t|_{t=t_{0}(u)}=0$, we get
\begin{equation*}
D_{u}r_{u}^{0}(t_{0}(u))(\phi )=t_{0}(u)\cdot
D_{w}R^{0}(w)|_{w=t_{0}(u)u}(\phi )=0,~~~\forall \phi \in C_{0}^{\infty
}(\Omega ).
\end{equation*}%
Now taking into account that the equality $\bar{\lambda}=\lambda _{0}(u)$
implies $E_{\bar{\lambda}}(u)=0,$ we obtain
\begin{equation*}
0=D_{w}R^{0}(w)|_{w=t_{0}(u)u}=\frac{1}{\int_{\Omega }|w|^{{\beta +1}}dx}%
\cdot D_{w}E_{\bar{\lambda}}(w)|_{w=t_{0}(u)u},
\end{equation*}%
which yields the proof.$\fin$

We shall need also the following Rayleigh's quotients:
\begin{align}
& R^{P}(u)=\frac{\frac{1}{2^{\ast }}\int_{\Omega }|\nabla u|^{2}\,\mathrm{d}%
x+\frac{1}{{\alpha +1}}\int_{\Omega }|u|^{{\alpha +1}}dx}{\frac{1}{{\beta +1}%
}\int_{\Omega }|u|^{{\beta +1}}dx},  \label{lamb11} \\
& R^{1}(u)=\frac{\int_{\Omega }|\nabla u|^{2}\,\mathrm{d}x+\int_{\Omega
}|u|^{{\alpha +1}}dx}{\int_{\Omega }|u|^{{\beta +1}}dx},~~u\neq 0.
\end{align}%
Notice that for any $u\neq 0$ and $\lambda \in \mathbb{R}$,
\begin{equation}
R^{P}(u)=\lambda \Leftrightarrow P_{\lambda }(u)=0~~\mbox{and}%
~~R^{1}(u)=\lambda \Leftrightarrow E_{\lambda }^{\prime }(u)=0.  \label{RPR1}
\end{equation}%
Let $u\neq 0$. Consider $r_{u}^{P}(t):=R^{P}(tu)$, $r_{u}^{1}(t):=R^{1}(tu)$%
, $t>0$. Then, arguing as above for $r_{u}^{0}(t),$ it can be shown that
each of these functions attains its global minimum at some point, $t_{P}(u)$
and $t_{1}(u)$, respectively. Moreover, it is easily seen that the following
equation
\begin{equation}
r_{u}^{P}(t)=r_{u}^{1}(t),~~~t>0,  \label{eqPD}
\end{equation}%
has a unique solution
\begin{equation}
t_{1P}(u)=\left( \frac{2^{\ast }(\beta -\alpha )}{(2^{\ast }-\beta
-1)(\alpha +1)}\frac{\int_{\Omega }|u|^{{\alpha +1}}dx}{\int_{\Omega
}|\nabla u|^{2}dx}\right) ^{\frac{1}{1-\alpha }}.  \label{P1}
\end{equation}%
Thus, we have {the next} \textit{nonlinear generalized Rayleigh quotient}
\begin{equation*}
\lambda _{1P}(u):=r_{u}^{P}(t_{1P}(u))=r_{u}^{1}(t_{1P}(u)).
\end{equation*}%
It is easily to seen that $\lambda _{1P}(u)=c_{1P}^{\alpha ,\beta }\lambda
(u)$, where
\begin{equation}
c_{1P}^{\alpha ,\beta }=\frac{(\beta +1)(2^{\ast }-\alpha +1)}{(\beta
-\alpha )2^{\ast }}\left( \frac{2^{\ast }(\beta -\alpha )}{(2^{\ast }-\beta
-1)(\alpha +1)}\right) ^{\frac{\beta -\alpha }{1-\alpha }}.  \label{cPD}
\end{equation}%
Notice that
\begin{equation}
P_{\lambda _{1P}(u)}(t_{1P}(u)u)=0,~~E_{\lambda
_{1P}(u)}^{\prime }(t_{1P}(u)u)=0,~~\forall u\neq 0.  \label{P1Propert}
\end{equation}%
Consider
\begin{equation}
\Lambda _{1P}=\inf_{u\neq 0}\lambda _{1P}(u).  \label{PPoh}
\end{equation}%
Using Sobolev's and H\"{o}lder's inequalities it can be shown (see, e.g.,
\cite{IlEg}) that
\begin{equation}
0<\Lambda _{1P}<+\infty .
\end{equation}

Moreover we have (see Figure 5):

%\begin{center}
%\FRAME{ftbpF}{5.2918in}{4.1978in}{0pt}{}{\Qlb{Fifure 2}}{figure 2.eps}{%
%\special{language "Scientific Word";type "GRAPHIC";maintain-aspect-ratio
%TRUE;display "USEDEF";valid_file "F";width 5.2918in;height 4.1978in;depth
%0pt;original-width 5.2356in;original-height 4.1476in;cropleft "0";croptop
%"1";cropright "1";cropbottom "0";filename 'Figure 2.eps';file-properties
%"XNPEU";}}
%\end{center}

\begin{claim}
\label{propEst} For any $u \neq 0$,

\begin{description}
\item[(i)] $r_{u}^{P}(t)>r_{u}^{1}(t)$ iff $t\in (0,t_{1P}(u))$ and $%
r_{u}^{P}(t)<r_{u}^{1}(t)$ iff $t\in (t_{1P}(u),+\infty )$;

\item[(ii)] $t_{1}(u)<t_{1P}(u)<t_{P}(u)$.
\end{description}
\end{claim}

\noindent {\em Proof}\quad Observe that $r_{u}^{P}(t)/r_{u}^{1}(t)\rightarrow
\frac{\beta +1}{{\alpha +1}}>1$ as $t\rightarrow 0$. Hence, from the
uniqueness of $t_{1P}(u)$ we obtain \textbf{(i)}.

\noindent By \eqref{RPR1} we have $E_{\lambda _{1P}(u)}^{\prime }(u)=0$.
Therefore Proposition \ref{pradd} implies $\frac{d}{dt}P_{\lambda
_{1P}(u)}(t_{1P}(u)u)<0$. Hence and since
\begin{equation*}
\frac{d}{dt}r_{u}^{P}(t)|_{t=t_{1P}(u)}=\frac{\beta +1}{\int |tu|^{{\beta +1}%
}dx}\cdot \frac{d}{dt}P_{\lambda _{1P}(u)}(tu)|_{t=t_{1P}(u)},
\end{equation*}%
we conclude that $\frac{d}{dt}r_{u}^{P}(t)|_{t=t_{1P}(u)}<0$. Now taking
into account that $t_{P}(u)$ is a point of global minimum of $r_{u}^{P}(t)$
we obtain that $t_{1P}(u)<t_{P}(u)$. To prove of $t_{1}(u)<t_{1P}(u)$, first
observe that
\begin{equation*}
\frac{d}{dt}r_{u}^{1}(t)|_{t=t_{1P}(u)}=\frac{1}{%
\int_{\Omega }|tu|^{{\beta +1}}dx}\cdot E_{\lambda_{1P}(u)}^{\prime \prime
}(tu)|_{t=t_{1P}(u)},
\end{equation*}%
and that by Lemma \ref{pro} the equalities $E_{\lambda _{1P}(u)}^{\prime
}(t_{1P}(u)u)=0$, $P_{\lambda _{1P}(u)}(t_{1P}(u)u)=0$ imply $E_{\lambda
_{1P}(u)}^{\prime \prime }(t_{1P}(u)u)>0$. Thus $\frac{d}{dt}%
r_{u}^{1}(t_{1P}(u))>0$ and the proof follows.$\fin$

\begin{cor}
\label{corPG}

\begin{description}
\item[(i)] If $\lambda<\Lambda_{1P}$ and $E_{\lambda}^{\prime }(u)=0$, then $%
P_\lambda(u)>0$.

\item[(ii) ] For any $\lambda>\Lambda_{1P}$, there exists $u \in
H_{0}^{1}\setminus 0 $ such that $E_{\lambda}^{\prime }(u)=0$ and $%
P_\lambda(u)\leq 0$
\end{description}
\end{cor}

\noindent {\em Proof}\quad \textbf{(i)}. ~Let $u \in H_{0}^{1}\setminus 0$. Assume $%
\lambda<\lambda_{1P}(u)$ such that $E_{\lambda}^{\prime }(u)=0$. Then in
view of \eqref{RPR1} we have $r^1_u(1)=\lambda<\lambda_{1P}(u)$. Thus
\textbf{(ii)}, Proposition \ref{propEst} yields $1\equiv t_{1}(u)< t_{1P}(u)$
and therefore by \textbf{(i)}, Proposition \ref{propEst} we have $%
r^P_u(1)>r^1_u(1)=\lambda$. Thus by \eqref{lamb1} we get $P_\lambda(u)>0$.

The proof of \textbf{(ii)} is similar to \textbf{(i)}.$\fin$

\begin{cor}
\label{proL} $\Lambda_{1P}<\Lambda_0$.
\end{cor}
\noindent {\em Proof}\quad
Suppose that $\Lambda _{0}<\Lambda _{1P}$. From
Proposition \ref{propL0} for any $\lambda \in (\Lambda _{0},\Lambda _{1P})$,
there exists $u \neq 0$ such that $E_{\lambda }(u)<0$ and $E_{\lambda
}^{\prime }(u)=0$. By Corollary \ref{corPG}, the equality $E_{\lambda
}^{\prime }(u)=0$ entails $P_{\lambda }(u)>0$. Hence by \eqref{PandE} we
have $E_{\lambda }(u)>P_{\lambda }(u)>0$, i.e., we get a contradiction. The
equality $\Lambda _{0}=\Lambda _{1P}$ is impossible since $c_{1P}^{\alpha
,\beta }\neq c_{0}^{\alpha ,\beta }.\fin$

\begin{cor}
\label{nonexist} Let $\Omega $ be a bounded star-shaped domain in $\mathbb{R}%
^{N}$ with $C^{2}$-manifold boundary $\partial \Omega $. Then for any $%
\lambda <\Lambda _{1P}$ equation $P(\alpha ,\beta ,\lambda )$ cannot have
weak solution.
\end{cor}

\noindent {\em Proof}\quad
Let $\lambda<\Lambda_{1P}$. Assume conversely that there
exists a weak solution $u$. By the regularity of solutions of elliptic
equations, $u \in C^1(\overline{\Omega})$. Then since $E_{\lambda}^{\prime
}(u)=0$ by Corollary \ref{corPG} we have $P_\lambda(u)>0$. However by
Corollary \ref{cor1}, any weak solution $u \in C^1(\overline{\Omega})$ of $%
P(\alpha ,\beta ,\lambda)$ satisfies $P_\lambda(u)\leq 0$. Thus we get a
contradiction. $\fin$

\medskip

%Observe, $r^P_u(t)>r^1_u(t)~~\mbox{if}~~t\in (0,t_{1P}(u)), ~\mbox{and}~~  r^1_u(t)<r^P_u(t)~~\mbox{if}~~t>t_{1P}(u)$. Using this  it is not hard to show that
%\begin{equation}\label{PDineq}
%\frac{d}{dt} r^1_u(t)<0,~~\forall t\in (0,t_{1P}(u)].
%\end{equation}

\section{Main constrained minimization problem}

Consider the constrained minimization problem:
\begin{equation}
\hat{E}_{\lambda }:=\min_{u\in M_{\lambda }}E_{\lambda }(u).  \label{min1}
\end{equation}%
where
\begin{equation*}
M_{\lambda }:=\{u\in H_{0}^{1}\setminus 0:~E_{\lambda }^{\prime
}(u)=0,~P_{\lambda }(u)\leq 0\}.
\end{equation*}%
Observe that any weak solution of $P(\alpha ,\beta ,\lambda )$ belongs to $%
M_{\lambda },$ such as it follows from Corollary \ref{cor1}. Hence if $\hat{E%
}_{\lambda }=E_{\lambda }(u_{\lambda })$, in \eqref{min1}, for some solution
$u_{\lambda }$ of $P(\alpha ,\beta ,\lambda )$, then $u_{\lambda }$ is a
ground state.

\begin{claim}
\label{PrL} $M_\lambda\neq \emptyset$ for any $\lambda> \Lambda_{1P}$.
\end{claim}

\noindent {\em Proof}\quad Let $\lambda >\Lambda _{1P}$. Consider the
function $\lambda _{1P}(\cdot ):H_{0}^{1}\setminus 0\rightarrow \mathbb{R}$.
By Proposition \ref{propDiff} this is a continuous functional. Hence there
is $u\in H_{0}^{1}\setminus 0$ such that $\Lambda _{1P}<\lambda
_{1P}(u)<\lambda $. Since by \eqref{P1Propert} we have $P_{\lambda
_{1P}(u)}(t_{1P}(u)u)=0$, $E_{\lambda _{1P}(u)}^{\prime }(t_{1P}(u)u)=0$, it
follows $P_{\lambda }(t_{1P}(u)u)<0$, $E_{\lambda }^{\prime }(t_{1P}(u)u)<0$%
. Hence there is $t_{\min }(u)>t_{1P}(u)$ such that $E_{\lambda }^{\prime
}(t_{\min }(u)u)=0$. In view that $P_{\lambda }^{\prime }(tu)=E_{\lambda
}^{\prime }(tu)-(2t/N)\int |\nabla u|^{2}$ for any $t>0$ we have $P_{\lambda
}^{\prime }(t_{\min }(u)u)<0$ which implies that $P_{\lambda }(t_{\min
}(u)u)<0$. Thus $t_{\min }(u)u\in M_{\lambda }.\fin$

%\subsection{Existence of the solution of (\ref{min1}).}

\begin{lem}
\label{le1e} For any $\lambda> \Lambda_{1P}$ there exists a minimizer $%
u_\lambda$ of problem (\ref{min1}), i.e., $E_\lambda(u_\lambda)=\hat{E}%
_\lambda$ and $u_\lambda\in M_\lambda$.
\end{lem}

\noindent {\em Proof}\quad Let $\lambda >\Lambda _{1P}$. Then $M_{\lambda }$
is bounded. Indeed, if $u\in M_{\lambda }$, %$$
%\frac{1}{2}\int |\nabla u|^{2}\,\mathrm{d}x+\frac{1%
%}{{\alpha +1}}\int_{\Omega }|u|^{{\alpha +1}}\,\mathrm{d}x-\lambda \frac{1}{{%
%\beta +1}}\int_{\Omega }|u|^{{\beta +1}}\,\mathrm{d}x\leq \frac{1}{N}\int |\nabla u|^{2}\,\mathrm{d}x
%$$
then
\begin{equation*}
\frac{1}{2^{\ast }}\int_{\Omega }|\nabla u|^{2}\,\mathrm{d}x+\frac{1}{{%
\alpha +1}}\int_{\Omega }|u|^{{\alpha +1}}\,\mathrm{d}x\leq \lambda \frac{1}{%
{\beta +1}}\int_{\Omega }|u|^{{\beta +1}}\,\mathrm{d}x\leq c\lambda \frac{1}{%
{\beta +1}}\Vert u\Vert _{1}^{\beta +1}
\end{equation*}%
From here $\Vert u\Vert _{1}\leq C<+\infty $, $\forall u\in M_{\lambda }$.
Now, if $(u_{m})$ is a minimizing sequence of (\ref{min1}), then it is
bounded and there exists a subsequence, denoting again, $(u_{m})$ which
converges $u_{m}\rightharpoonup u_{0}$ weakly in $H_{0}^{1}$ and strongly $%
u_{m}\rightarrow u_{0}$ in $L^{q}$, $1<q<2^{\ast }$. We claim that $%
u_{m}\rightarrow u_{0}$ strongly in $H_{0}^{1}$. If not, $\Vert u_{0}\Vert
_{1}<\liminf_{m\rightarrow \infty }\Vert u_{m}\Vert _{1}$ and this implies
\begin{align*}
\int_{\Omega }|\nabla u_{0}|^{2}\,\mathrm{d}x+& \int_{\Omega }|u_{0}|^{{%
\alpha +1}}\,\mathrm{d}x-\lambda \int_{\Omega }|u_{0}|^{{\beta +1}}\,\mathrm{%
d}x< \\
& \liminf_{m\rightarrow \infty }\left( \int_{\Omega }|\nabla u_{m}|^{2}\,%
\mathrm{d}x+\int_{\Omega }|u_{m}|^{{\alpha +1}}\,\mathrm{d}x-\lambda
\int_{\Omega }|u_{m}|^{{\beta +1}}\,\mathrm{d}x\right) =0
\end{align*}%
since $E_{\lambda }^{\prime }(u_{m})=0$, $m=1,2,...$. Hence $u_{0}\neq 0$
and $E_{\lambda }^{\prime }(u_{0})=0$. Then there exists $\gamma >1$ such
that $E_{\lambda }^{\prime }(\gamma u_{0})=0$ and $E_{\lambda }(\gamma
u_{0})<E_{\lambda }(u_{0})<\hat{E}_{\lambda }$. By Proposition \ref{pradd}, $%
E_{\lambda }^{\prime }(\gamma u_{0})=0$ implies $P_{\lambda }^{\prime
}(\gamma u_{0})<0$. From this and since
\begin{equation*}
P_{\lambda }(u_{0})<\liminf_{m\rightarrow \infty }P_{\lambda }(u_{m})\leq 0,
\end{equation*}%
we conclude that $P_{\lambda }(\gamma u_{0})<0$. Thus $\gamma u_{0}\in
M_{\lambda }$ and $E_{\lambda }(\gamma u_{0})<\hat{E}_{\lambda }$, which is
a contradiction.$\fin$

\subsection{Existence of a flat or compact support ground state $u_{\protect%
\lambda ^{\ast }}$}

Let $\lambda >\Lambda _{1P}$, then by Lemma \ref{le1e} there exists
a minimizer $u_{\lambda }$ of (\ref{min1}). Notice since $\min
\{\alpha ,\beta \}>0$, $E_{\lambda }(u)$ and $E_{\lambda }^{\prime
}(u)$, $P_{\lambda }(u)$ are $C^{1}$-functionals on
$H_{0}^{1}(\Omega )$ . Hence we may apply Lagrange multipliers rule
(see, e.g., page 417 of \cite{Zeidler}) and thereby there exist
Lagrange multipliers $\mu _{0}$, $\mu _{1}$ $\mu _{2}$ such that
$|\mu _{0}|+|\mu _{1}|+|\mu _{2}|\neq 0$, $\mu _{2}\geq 0$ and
\begin{eqnarray}
&&\mu _{0}D_{u}E_{\lambda }(u_{\lambda })+\mu _{1}D_{u}E_{\lambda }^{\prime
}(u_{\lambda })+\mu _{2}D_{u}P_{\lambda }(u_{\lambda })=0,  \label{eq2} \\
&&\mu _{2}P_{\lambda }(u_{\lambda })=0.  \label{eq22}
\end{eqnarray}

\begin{claim}
\label{Lag} Assume $(\alpha ,\beta )\in \mathcal{E}_{s}(N)$. Let $\lambda
>\Lambda _{1P}$ and $u_{\lambda }\in H_{0}^{1}$ be a minimizer in (\ref{min1}%
) such that $P_{\lambda }(u_{\lambda })<0$. Then $u_{\lambda }$ is a weak
solution to $P(\alpha ,\beta ,\lambda )$.
\end{claim}

\noindent {\em Proof}\quad Since $P_{\lambda }(u_{\lambda })<0$, equality %
\eqref{eq22} implies $\mu _{2}=0$. Moreover, since $(\alpha ,\beta
)\in \mathcal{E}_{s}(N)$, ($ii$), Lemma \ref{pro} implies that
$E_{\lambda }^{\prime \prime }(u_{\lambda })>0$. Testing \eqref{eq2}
by $u_{\lambda }$ we get $0=\mu _{0}E_{\lambda }^{\prime
}(u_{\lambda })=\mu _{1}E_{\lambda }^{\prime \prime }(u_{\lambda
})$. But $E_{\lambda }^{\prime \prime }(u_{\lambda })\neq 0$ and
therefore $\mu _{1}=0$. Thus, $D_{u}E_{\lambda }(u_{\lambda })=0$,
that is $u_{\lambda }$ weakly satisfies $P(\alpha ,\beta ,\lambda
)$. This completes the proof.$\fin$

\noindent Introduce
\begin{equation}
Z:=\{\lambda \in (\Lambda _{1P},+\infty ):~~P_{\lambda }(u_{\lambda
})<0,~~u_{\lambda }\in M_{\lambda }~~\mbox{s.t.}~E_{\lambda }(u_{\lambda })=%
\hat{E}_{\lambda }\}.  \label{Pr1}
\end{equation}

\begin{claim}
$Z$ is a non-empty open subset of $(\Lambda_{1P}, +\infty)$.
\end{claim}

\noindent {\em Proof}\quad Notice that by Lemma \ref{le1e}, for any $\lambda
>\Lambda _{1P}$ there exists $u_{\lambda }\in M_{\lambda }$ such that $%
E_{\lambda }(u_{\lambda })=\hat{E}_{\lambda }$. To prove that $Z\neq
\emptyset $, we show that $[\Lambda _{0},+\infty )\subset Z$. Take $\lambda
\geq \Lambda _{0}$. Then in view of (ii), Proposition \ref{propL0} we have $%
\hat{E}_{\lambda }\leq 0$. Thus $E_{\lambda }(u_{\lambda })\leq 0$, for any $%
u_{\lambda }\in M_{\lambda }$ such that $E_{\lambda }(u_{\lambda })=\hat{E}%
_{\lambda }$. In view of \eqref{PandE} we have $E_{\lambda }(u_{\lambda
})>P_{\lambda }(u_{\lambda })$ and therefore $P_{\lambda }(u_{\lambda })<0$,
$\forall u_{\lambda }\in M_{\lambda }$~~ such that $E_{\lambda }(u_{\lambda
})=\hat{E}_{\lambda }$. Thus $\lambda \in Z$.

\noindent Let us show that $Z$ is an open subset of $(\Lambda _{1P},+\infty )
$. Notice that if $Z=(\Lambda _{1P},+\infty )$, then $Z$ is an open subset
of $(\Lambda _{1P},+\infty )$ by the definition.

\noindent Assume $Z\neq (\Lambda _{1P},+\infty )$. Let $\lambda \in Z$.
Suppose, contrary to our claim, that there is a sequence $(\lambda
_{m})\subset (\Lambda _{1P},+\infty )\setminus Z$ such that $\lambda
_{m}\rightarrow \lambda $ as $m\rightarrow \infty $. Then there is a
sequences of solutions $(u_{\lambda _{m}})$ of (\ref{min1}) such that $%
P_{\lambda _{m}}(u_{\lambda _{m}})=0$. Then by Lemma \ref{app} (see Appendix
I), there exists a minimizer $u_{\lambda }$ of (\ref{min1}) and a
subsequence, still denoted by $(u_{\lambda _{m}})$, such that $u_{\lambda
_{m}}\rightarrow u_{\lambda }$ strongly in $H_{0}^{1}$ as $m\rightarrow
+\infty $. However, then $P_{\lambda }(u_{\lambda })=0$, which contradicts
the assumption $\lambda \in Z.\fin$

Set
\begin{equation*}
\lambda^*:=\inf Z.
\end{equation*}

\begin{lem}
\label{Pr2} There exists a minimizer $u_{\lambda ^{\ast }}$ of (\ref{min1})
which is a flat or a compact support non-negative ground state of $P(\alpha
,\beta ,\lambda ^{\ast })$. Furthermore, $\Lambda _{1P}<\lambda ^{\ast }$
and there exists a set of \textquotedblleft usual\textquotedblright\
non-negative ground states $(u_{\lambda _{n}})_{n=1}^{\infty }$ of $P(\alpha
,\beta ,\lambda _{n})$, with $\lambda _{n}\downarrow \lambda ^{\ast }$ as $%
n\rightarrow \infty $, such that $u_{\lambda _{n}}\rightarrow u_{\lambda
^{\ast }}$ strongly in $H_{0}^{1}$ as $n\rightarrow \infty $.
\end{lem}

\noindent {\em Proof}\quad Since $Z$ is an open set, we can find a sequence $%
\lambda _{n}\in Z$, $n=1,2,...$ such that $\lambda _{n}\rightarrow \lambda
^{\ast }$ as $n\rightarrow \infty $. By the definition of $Z$ for any $%
n=1,2,...$ we can find a minimizer $u_{\lambda _{n}}$ of (\ref{min1}) such
that $P_{\lambda _{n}}(u_{\lambda _{n}})<0$. Then Proposition \ref{Lag}
yields that $u_{\lambda _{n}}$ weakly satisfies $P(\alpha ,\beta ,\lambda
_{n})$, $n=1,2,...$. Moreover by Corollary \ref{cor1}, $u_{\lambda _{n}}$ is
a \textquotedblleft usual\textquotedblright\ weak solution of $P(\alpha
,\beta ,\lambda _{n})$, $n=1,2,...$. Since $E_{\lambda }(|u|)=E_{\lambda }(u)
$, $E_{\lambda }^{\prime }(|u|)=E_{\lambda }^{\prime }(u)=0$, $P_{\lambda
}(|u|)=P_{\lambda }(u)$ for any $u\in H_{0}^{1}$ we may assume that $%
u_{\lambda _{n}}\geq 0$, $n=1,2,...$. Furthermore, since $\hat{E}_{\lambda
_{n}}={E}_{\lambda _{n}}(u_{\lambda _{n}})$, $u_{\lambda _{n}}$ is a ground
state of $P(\alpha ,\beta ,\lambda _{n})$, $n=1,2,...$. Thus we have a set
of \textquotedblleft usual\textquotedblright\ non-negative ground states $%
(u_{\lambda _{n}})_{n=1}^{\infty }$ of $P(\alpha ,\beta ,\lambda _{n})$, $%
n=1,2,...$.

\noindent By Lemma \ref{app} (see Appendix I), there exists a minimizer $%
u_{\lambda ^{\ast }}$ of (\ref{min1}) and the subsequence, still denoted by $%
(u_{\lambda _{n}})$, such that $u_{\lambda _{n}}\rightarrow u_{\lambda
^{\ast }}$ strongly in $H_{0}^{1}$ as $\lambda _{n}\rightarrow \lambda
^{\ast }$. This implies that $u_{\lambda ^{\ast }}$ is a non-negative
solutions of $P(\alpha ,\beta ,\lambda )$ and $P_{\lambda ^{\ast
}}(u_{\lambda ^{\ast }})\leq 0$. Furthermore, since $u_{\lambda ^{\ast }}$
is a minimizer of (\ref{min1}), it is a ground state of $P(\alpha ,\beta
,\lambda ^{\ast })$.

\noindent Let us show that $\Lambda _{1P}<\lambda ^{\ast }$. To obtain a
contradiction suppose, that $\Lambda _{1P}=\lambda ^{\ast }$. Then $\Lambda
_{1P}=\lambda _{1P}(u_{\lambda ^{\ast }})$ and $u_{\lambda ^{\ast }}$ is a
minimizer of \eqref{PPoh}. Since $\lambda _{1P}(u)=c^{\alpha ,\beta }\lambda
_{0}(u)$, where $c^{\alpha ,\beta }=c_{1P}^{\alpha ,\beta }/c_{0}^{\alpha
,\beta }$, $u_{\lambda ^{\ast }}$ is also a critical point of $\lambda
_{0}(u)$ with value $\Lambda _{0}$. Then by Proposition \ref{CRRay}, $%
u_{\lambda ^{\ast }}$ satisfies $P(\alpha ,\beta ,\Lambda _{0})$. However,
by the construction $u_{\lambda ^{\ast }}$ satisfies $P(\alpha ,\beta
,\lambda ^{\ast })$. Notice that by Corollary \ref{proL}, $\Lambda
_{0}>\Lambda _{1P}=\lambda ^{\ast }$. Thus we get a contradiction.

\noindent Observe that $P_{\lambda ^{\ast }}(u_{\lambda ^{\ast }})=0$.
Indeed, if $P_{\lambda ^{\ast }}(u_{\lambda ^{\ast }})<0$, then $\lambda
^{\ast }\in Z$. But this is impossible since $Z$ is an open subset of $%
(\Lambda _{1P},+\infty )$.

\noindent A global (up to the boundary) regularity result (see \cite%
{Lieberman}) yields that $u_{\lambda}\in C^{1,\beta
}(\overline{\Omega })$, $\lambda \in \lbrack \lambda ^{\ast },+\infty )$ for
some $\beta \in (0,1)$. Thus we may apply Corollary \ref{cor1} which yields
that $u_{\lambda ^{\ast }}$ is flat or compactly supported in $\Omega .\fin$

\bigskip

%\begin{rem}
%From the above it follows the existence of non-negative ground state $%
%u_\lambda$ of $P(\alpha ,\beta ,\lambda)$ for any $\lambda \in Z$. Thus to
%conclude the proof Theorem \ref{Th1} it is remain to show that $%
%Z=(\lambda^*, +\infty)$.
%\end{rem}

%Let us show that for any $\lambda>\lambda^*$ the ground state $u_\lambda$  can not be compact supported. First, notice that $\hat{E}_\lambda> \hat{E}_\mu$ for any $\lambda<\mu$, $\lambda,\mu \in (\Lambda_{1P},+\infty)$. Indeed, $\hat{E}_\lambda=E_\lambda(u_\lambda)>E_\mu(u_\lambda)$. Furthermore, since $E_\mu'(u_\lambda)<0$, $P_\mu(u_\lambda)<0$, there is $t_\mu(u_\lambda)>1$ such that  $E_\mu'(t_\mu(u_\lambda)u_\lambda)=0$, $P_\mu(t_\mu(u_\lambda)u_\lambda)<0$ and $E_\mu(u_\lambda)> E_\mu(t_\mu(u_\lambda)u_\lambda)$. Hence $\hat{E}_\lambda>E_\mu(t_\mu(u_\lambda)u_\lambda)\geq  \hat{E}_\mu$.

\section{On the radially symmetric property}

We need the following result that has been proved in \cite{Kaper1, Kaper2,
Serrin-Zou}.

\begin{lem}
\label{lem:3} Assume $0<\alpha <\beta<1$. Let $u$ be a non-negative $C^1$
distribution solution of
\begin{equation}  \label{Eqw}
-\Delta u+u^{\alpha}=u^{\beta}~~~\mbox{in}~~\mathbb{R}^N
\tag*{$Eq(\alpha
,\beta ,1)$}
\end{equation}
with connected support. Then the support of $u$ is a ball and $u$ is
radially symmetric about the center .

Furthermore, equation $Eq(\alpha ,\beta ,1)$ admits at most one radially
symmetric compact support solution.
\end{lem}

We denote by $R^{\ast }$ the radius of the supporting ball $B_{R^{\ast }}$
of the unique (up to translation in $\mathbb{R}^{N}$) compact support
solution of $Eq(\alpha ,\beta ,1)$, i.e., it is the unique flat solution of $%
P(\alpha ,\beta ,1)$ for $\Omega =B_{R^{\ast }}$.

It is easy to see, from Lemma \ref{lem:3}, that the function $u_{\lambda
}^{\ast }(x):=\sigma ^{-\frac{2}{1-\alpha }}\cdot u^{\ast }(x/\sigma )$ is
the unique flat solution of
$P(\alpha ,\beta ,\lambda)$  with
$\lambda =\sigma ^{-\frac{2(\beta -\alpha )}{1-\alpha }}$ and
$\Omega=B_{\sigma R^{\ast }}$.

\begin{claim}
Assume $u_{\lambda }\in C^{1}(\overline{\Omega })$ is a non-negative ground
state of $P(\alpha ,\beta ,\lambda )$ which has compact support in $\Omega $%
. Then $u_{\lambda }$ is radially symmetric about some origin $0\in \Omega $%
, and supp$(u_{\lambda })$=$\overline{B_{R(\Omega )}}$ is a inscribed ball
in $\overline{\Omega }$.
\end{claim}

\noindent {\em Proof}\quad Observe that any compact support function $%
u_{\lambda }$ from $C^{1}(\overline{\Omega })$ can be extended to $\mathbb{R}%
^{N}$ as
\begin{equation}
\left\{
\begin{array}{ll}
\tilde{u}_{\lambda }=u_{\lambda } & \mbox{in}~\Omega , \\
\tilde{u}_{\lambda }=0 & \mbox{in}~\mathbb{R}^{N}\setminus \Omega .%
\end{array}%
\right.   \label{expan}
\end{equation}%
Then $\tilde{u}_{\lambda }\in C^{1}(\mathbb{R}^{N})$ is a distribution
solution of $P(\alpha ,\beta ,\lambda )$ on $\mathbb{R}^{N}$. Since $%
u_{\lambda }$ is a ground state, it is not hard to show that $u_{\lambda }$
has a connected support. Thus by Lemma \ref{lem:3}, $\tilde{u}_{\lambda }$
is a radially symmetric function with respect to the centre of some ball $%
B_{R^{\lambda }}$ with a radius $R^{\lambda }>0,$ so that supp($u_{\lambda
})=\overline{B}_{R^{\lambda }}$.

\noindent Let us show that $B_{R^{\lambda }}$ is an inscribed ball in $%
\Omega $. Consider $B_{\sigma R^{\lambda }}:=\{x\in \mathbb{R}^{N}:~x/\sigma
\in B_{R^{\lambda }}\}$ where $\sigma >0$. Notice that $B_{\sigma R^{\lambda
}}\subset \Omega $ if $\sigma \leq 1$. Suppose, contrary to our claim, that
there is $\sigma _{0}>1$ such that $B_{\sigma R^{\lambda }}\subset \Omega $
for any $\sigma \in (1,\sigma _{0})$. Let $\sigma \in (1,\sigma _{0})$.
Introduce $u_{\lambda }^{\sigma }(x)=u_{\lambda }(x/\sigma ),~x\in B_{\sigma
R^{\lambda }}$ and set $u_{\lambda }^{\sigma }(x)=0$ in $\Omega \setminus
B_{\sigma R^{\lambda }}$. Observe that
\begin{equation*}
E_{\lambda }(u_{\lambda }^{\sigma })=\frac{\sigma ^{N-2}}{2}\int_{\Omega
}|\nabla u_{\lambda }|^{2}\,\mathrm{d}x-\sigma ^{N}(\frac{\lambda }{{\beta +1%
}}\int_{\Omega }|u_{\lambda }|^{{\beta +1}}\,\mathrm{d}x-\frac{1}{{\alpha +1}%
}\int_{\Omega }|u_{\lambda }|^{{\alpha +1}}\,\mathrm{d}x).
\end{equation*}%
From this $dE_{\lambda }(u_{\lambda }^{\sigma })/d\sigma |_{\sigma
=1}=P_{\lambda }(u_{\lambda })=0$, and thus $\sigma =1$ is a maximizing
point of the function $\psi _{u}(\sigma ):=E_{\lambda }(u_{\lambda }^{\sigma
})$. Then $E_{\lambda }(u_{\lambda }^{\sigma })<E_{\lambda }(u_{\lambda })=%
\hat{E}_{\lambda }$ and $P_{\lambda }(u_{\lambda }^{\sigma })<0$ for $\sigma
\in (1,\sigma _{0})$. From this it follows that for $\sigma $ sufficiently
close to $1$ we  have $E_{\lambda }(t_{\min }(u_{\lambda }^{\sigma
})u_{\lambda }^{\sigma })<E_{\lambda }(u_{\lambda })=\hat{E}_{\lambda }$ and
$P_{\lambda }(t_{\min }(u_{\lambda }^{\sigma })u_{\lambda }^{\sigma })<0$, $%
E_{\lambda }^{\prime }(t_{\min }(u_{\lambda }^{\sigma })u_{\lambda }^{\sigma
})=0,$ which is a contradiction. $\fin$

\noindent From this and Lemma \ref{Pr2} we have

\begin{cor}
\label{corFin1} $u_{\lambda ^{\ast }}$ is radially symmetric about some
point of $\Omega $, and supp$(u_{\lambda ^{\ast }})$=$\overline{B}_{R(\Omega
)}$ is an inscribed ball in $\Omega $.
\end{cor}

Furthermore, we have

\begin{cor}
\label{corFin2} For any $\lambda > \lambda^*$, problem $P(\alpha ,\beta
,\lambda)$ has no non-negative ground state with compact support.
\end{cor}

\noindent {\em Proof}\quad Suppose, conversely that there exists ${\lambda_a}%
>\lambda^*$ and a ground state $u_{\lambda_a}$ of $P(\alpha ,\beta
,\lambda_a)$ such that $u_{\lambda_a}$ has a compact support. Then arguing
as above one may infer that $u_{{\lambda_a}}$ is a radially symmetric
function with respect to a centre of inscribed ball $B_{R(\Omega)}$ in $%
\Omega$ so that supp($u_{\lambda_a})=\overline{B}_{R(\Omega)}$. Consider $%
u_{\lambda^*}^\sigma(x)=u_{\lambda^*}(x/\sigma)$ with $\sigma=(\lambda^*/{%
\lambda_a})^{(1-\alpha)/2(\beta-\alpha)}$. Then $u_{\lambda^*}^\sigma$ is
compactly supported non-negative weak solution of $P(\alpha ,\beta
,\lambda_a)$. By the uniqueness of radial compact support solution of $%
P(\alpha ,\beta ,\lambda_a)$ (see Lemma \ref{lem:3}) this is possible only
if $u_{\lambda^*}^\sigma=u_{\lambda_a}$. However supp($u_{\lambda^*}^\sigma)=%
\overline{B}_{\sigma R(\Omega)}$ whereas supp($u_{\lambda_a})=\overline{B}%
_{R(\Omega)}$ and $\sigma<1$. Thus we get a contradiction.$\fin$

\begin{cor}
\label{corFin3} $Z=(\lambda^*, +\infty)$.
\end{cor}

\noindent {\em Proof}\quad Suppose, contrary to our claim, that there is an
another limit point $\lambda _{b}$ of $Z$ such that $\lambda _{b}\in
(\lambda ^{\ast },+\infty )\setminus Z$. Then arguing similarly to the proof
of Lemma \ref{Pr2} one may conclude that there exists a compactly supported
non-negative ground state $u_{\lambda _{b}}$ of $P(\alpha ,\beta ,\lambda
_{b})$. However $\lambda _{b}>\lambda ^{\ast }$ and therefore by Corollary %
\ref{corFin2} this is impossible.$\fin$

%Moreover, since $(1^o)$, Lemma \ref{pro}, any compact supported weak solution $u$ of $P(\alpha ,\beta ,\lambda)$ satisfies $E_\lambda(u)>0$,
%
%
%
%Let us show that for $\lambda\in (\Lambda_{1P}, \lambda^*)$ any minimizer $u_\lambda \in M_\lambda$  can not satisfy $P(\alpha ,\beta ,\lambda)$. Suppose this is false, that is there exists $\lambda\in (\Lambda_{1P}, \lambda^*)$ such that $u_\lambda$ satisfies $P(\alpha ,\beta ,\lambda)$. Then since $\lambda^*=\inf Z$, this is possible only if $P_\lambda(u_\lambda)=0$.

\section{Proofs of Theorems}

\subsection{Proof of Theorem \protect\ref{Th1}}

For $\lambda=\lambda^*$, the existence of non-negative ground state $%
u_{\lambda^*}$ of $P(\alpha ,\beta ,\lambda^*)$ follows from Lemma \ref{Pr2}%
. Since $Z=(\lambda^*, +\infty)$, we see that for $\lambda>\lambda^*$, any
minimizer $u_\lambda $ of (\ref{min1}) satisfies $P_\lambda(u_\lambda)< 0 $.
From this by Proposition \ref{Lag} we derive that $u_\lambda $ is a weak
solution of $P(\alpha ,\beta ,\lambda)$. Moreover, since $\hat{E}_\lambda={E}%
_\lambda(u_\lambda)$, $u_\lambda$ is a ground state of $P(\alpha ,\beta
,\lambda )$ for all $\lambda \in (\lambda^*, +\infty)$. By the same
arguments as in the proof of Lemma \ref{Pr2} we may assume that $%
u_{\lambda}\geq 0$ in $\Omega$ for all $\lambda> \lambda^*$. In view of
Lemma \ref{pro} we have $E^{\prime \prime }_\lambda(u_{\lambda})>0$, and by
global (up to the boundary) regularity result for elliptic equations we have
$u_\lambda \in C^{1,\gamma}(\overline{\Omega})$ for some $\gamma\in (0,1)$.

Let us prove that for $\lambda < \lambda^*$, problem $P(\alpha ,\beta
,\lambda)$ has no weak solution $u \in H^1_0(\Omega)$. Observe that any weak
solution of $P(\alpha ,\beta ,\lambda)$ (if it exists) by global (up to the
boundary) regularity result for elliptic equations belongs to $C^1(\overline{%
\Omega})$. Notice that by Corollary \ref{nonexist} for any $%
\lambda<\Lambda_{1P}$ equation $P(\alpha ,\beta ,\lambda)$ has no weak
solution $u \in C^1(\overline{\Omega}) $. Thus since by Lemma \ref{Pr2}, $%
\Lambda_{1P}<\lambda^*$ it remains to prove nonexistence of weak solutions
in the case $\lambda \in [\Lambda_{1P}, \lambda^*)$.

Let $\lambda \in [\Lambda_{1P}, \lambda^*)$. Suppose, contrary to our claim,
that there exists a weak solution $u_\lambda \in C^1(\overline{\Omega})$ of $%
P(\alpha ,\beta ,\lambda)$. Then $E^{\prime}(u_\lambda)=0$ and by Corollary %
\ref{cor1} we have $P_\lambda(u_\lambda)\leq 0$. Hence $u_\lambda \in
M_\lambda$.

Let us show that then there exists a ground state of $P(\alpha ,\beta
,\lambda )$ which belongs to $C^{1}(\overline{\Omega })$. Notice that if $%
u_{\lambda }$ is a unique solution of $P(\alpha ,\beta ,\lambda )$ then it
is a ground state. Assume there exists a set of such solutions $\tilde{M}%
_{\lambda }$ of $P(\alpha ,\beta ,\lambda )$. Notice that $\tilde{M}%
_{\lambda }\subset M_{\lambda }$. Consider
\begin{equation}
\tilde{E}_{\lambda }:=\min_{u\in \tilde{M}_{\lambda }}E_{\lambda }(u).
\label{min1Til}
\end{equation}%
%
%
%
%
%
%
%
%
%
%
%
%
%
%Observe that in view of Lemma \ref{le1e}, for $\lambda\geq \Lambda_{1P}$, there holds $%
%\tilde{E}_\lambda\geq \hat{E}_\lambda>-\infty$.
Let $(u_{m})$ be a minimizing sequence of (\ref{min1Til}), i.e.,
\begin{equation}
E_{\lambda }(u_{m})\rightarrow \tilde{E}_{\lambda }~~\mbox{as}~~n\rightarrow
\infty ~~\mbox{and}~u_{m}\in \tilde{M}_{\lambda },~n=1,2,...  \label{maxsTil}
\end{equation}%
Using the same arguments as in the proof of Lemma \ref{le1e} we may conclude
that there exists a nonzero limit point $\tilde{u}_{0}$ such that (up to
subsequence) $u_{m}\rightarrow \tilde{u}_{0}$ converges weakly in $H_{0}^{1}$
and strongly in $L_{q}$ for $1<q<2^{\ast }$. Then we have
\begin{equation}
E_{\lambda }(\tilde{u}_{0})\leq \tilde{E}_{\lambda }  \label{lower}
\end{equation}%
and
\begin{equation*}
0=D_{u}E_{\lambda }(u_{m})(\psi )\rightarrow D_{u}E_{\lambda }(\tilde{u}%
_{0})(\psi )~~~\forall \psi \in C_{0}^{\infty }(\Omega ).
\end{equation*}%
Thus $\tilde{u}_{0}$ is a nonzero weak solution of $P(\alpha ,\beta ,\lambda
)$. Moreover by global (up to the boundary) regularity result for elliptic
equations we have $\tilde{u}_{0}\in C^{1,\gamma }(\overline{\Omega })$ for
some $\gamma \in (0,1)$. Thus $\tilde{u}_{0}\in \tilde{M}_{\lambda }$ and by %
\eqref{lower} we conclude that $E_{\lambda }(\tilde{u}_{0})=\tilde{E}%
_{\lambda }$. This implies that $\tilde{u}_{0}$ is a ground state of $%
P(\alpha ,\beta ,\lambda )$ belonging to $C^{1}(\overline{\Omega })$.

Thus we have proved that there exists a ground state $u_{\lambda }$ of $%
P(\alpha ,\beta ,\lambda )$ which belongs to $C^{1}(\overline{\Omega })$.
Then there are two possibilities $P_{\lambda }(u_{\lambda })<0$ or $%
P_{\lambda }(u_{\lambda })=0$. In the first case, we get that $\lambda \in Z$%
. But in view of Corollary \ref{corFin3} this is a contradiction. In the
second case, Corollary \ref{cor1} implies that $u_{\lambda }$ has a compact
support in $\Omega $. However the same arguments as in the proof of
Corollary \ref{corFin2} show that for $\lambda \neq \lambda ^{\ast }$ this
is impossible.

This concludes the proof of Theorem \ref{Th1}.

\subsection{Proof of Theorem \protect\ref{Th2}}

The existence of a non-negative ground state $u_{\lambda ^{\ast }}$ with
compact support follows from Lemma \ref{Pr2} . By Corollary \ref{corFin1}, $%
u_{\lambda ^{\ast }}$ is radially symmetric about some point of $\Omega $,
and supp$(u)$=$\overline{B}_{R(\Omega)}$ is an
inscribed ball in $\Omega $.

In view of Corollary \ref{corFin2}, for all $\lambda > \lambda^*$, any
ground state $u_\lambda$ of $P(\alpha ,\beta ,\lambda)$ is a usual solution.

\subsection{Proof of Theorem \protect\ref{ThmCor1}}

We shall only prove the theorem, as an example, for the case $m=2$, i.e.,
when $\Omega $ is a domain of Strictly Starshaped Class $2$.

Let $\lambda ^{\ast }>0$ be a limit value obtained in Theorem \ref{Th1}. By
Lemma \ref{Pr2} there exists a compactly supported ground state $u_{\lambda
^{\ast }}^{1}$ of $P(\alpha ,\beta ,\lambda )$ and there exists a set of
usual non-negative ground states $(u_{\lambda _{n}}^{1})_{n=1}^{\infty }$,%
 $\lambda_n>\lambda^{\ast }$, $n=1,2,...$ such that $%
u_{\lambda_n}^{1}\rightarrow u_{\lambda ^{\ast }}^{1}$ strongly in $%
H_{0}^{1} $ as $n\rightarrow \infty $. By Corollary \ref{corFin1}, $%
u_{\lambda^*}$ is radially symmetric about some origin $0\in \Omega$, and
supp$(u)$=$\overline{B}_{R(\Omega)}$ is an inscribed
ball in $\Omega$. By the assumptions $\Omega$ contains exactly $2$ inscribed
balls of radio $R(\Omega) $

Set $u_{\lambda ^{\ast }}^{2}(x):=u_{\lambda ^{\ast }}^{1}(R_{H}x)$, $%
u_{\lambda _{n}}^{2}(x):=u_{\lambda _{n}}^{1}(R_{H}x)$, $x\in \Omega $, $%
n=1,2,...$, where $R_{H}:\mathbb{R}^{N}\rightarrow \mathbb{R}^{N}$ is the
reflection map. By Theorem \ref{Th1}, the support of $u_{\lambda ^{\ast
}}^{1}$ coincides with one of the balls $B^{1}$ or $B^{2}$. Assume supp($%
u_{\lambda ^{\ast }}^{1})=B^{1}$. Then since $R_{H}B_{1}=B_{2}$ for some
hyperplane $H$, we have supp($u_{\lambda ^{\ast }}^{2})=B^{2}$ and thus $%
u_{\lambda ^{\ast }}^{2}\neq u_{\lambda ^{\ast }}^{1}$. Since $u_{\lambda
_{n}}^{2}\rightarrow u_{\lambda ^{\ast }}^{2}$ strongly in $H_{0}^{1}$ as $%
n\rightarrow \infty $, it follows that $u_{\lambda _{n}}^{1}\neq u_{\lambda
_{n}}^{2}$ for sufficiently large $n$.

\bigskip

\section{On the free boundary for the parabolic problem}

We consider now the associate parabolic problem

\begin{equation}
PP(\alpha ,\beta ,\lambda ,v_{0})\quad \left\{
\begin{array}{ll}
v_{t}-\Delta v+|v|^{\alpha -1}v=\lambda |v|^{\beta -1}v & \text{in }%
(0,+\infty )\times \Omega \\
v=0 & \text{on }(0,+\infty )\times \partial \Omega \\
v(0,x)=v_{0}(x) & \text{on }\Omega .%
\end{array}%
\right.
\end{equation}%
For the basic theory for this problem, always under the structural
assumption $0<\alpha <\beta <1$, we send the reader to \cite{CDE}
and  \cite{DIH}. In
particular, we know that for any $v_{0}\in \mathrm{L}^{\infty }(\Omega ),$ $%
v_{0}\geq 0$ there exists a nonnegative weak solution $v\in \mathcal{C}%
([0,+\infty ),\mathrm{L}^{2}(\Omega ))\cap $ $L^{\infty }((0,+\infty )\times
\Omega )$ of $PP(\alpha ,\beta ,\lambda ,v_{0})$. This solution is unique if
$v_{0}$ is non-degenerate near its free boundary.

Our main goal in this Section is to give an idea of the time evolution of
the support of the solution. We recall that, as $t\rightarrow +\infty $, the
support of $v(t,.)$ must converge to a ball of $\mathbb{R}^{N}$, in the case
$\lambda =\lambda ^{\ast }$, or to the whole domain $\overline{\Omega }$, if
$\lambda >\lambda ^{\ast }$ (since the shape of the support of the
associated stationary solutions was given in Theorem 1.1).

Our first result concerns the special case of $v_{0}=u_{\lambda ^{\ast }}$
(i.e. with support in the ball of $\mathbb{R}^{N}$ of radio $R(\Omega )$)
and $\lambda >\lambda ^{\ast }$. It is clear that any stationary solution $%
u_{\lambda ^{\ast }}$ is a subsolution to the problem $PP(\alpha ,\beta
,\lambda ,v_{0})$. Indeed,
\begin{equation*}
(u_{\lambda ^{\ast }})_{t}-\Delta u_{\lambda ^{\ast }}+|u_{\lambda ^{\ast
}}|^{\alpha -1}u_{\lambda ^{\ast }}=\lambda ^{\ast }|u_{\lambda ^{\ast
}}|^{\beta -1}u_{\lambda ^{\ast }}<\lambda |u_{\lambda ^{\ast }}|^{\beta
-1}u_{\lambda ^{\ast }}.
\end{equation*}%
So, if $u_{\lambda ^{\ast }}$ is nondegenerate near its free boundary, we
get that $u_{\lambda ^{\ast }}(x)\leq v(t,x)$ for any $t>0$ and a.e. $x\in
\Omega .$ As a matter of fact, it is easy to prove that under these
assumptions $v_{t}\geq 0$ a.e. $(0,+\infty )\times \Omega $. Thus, a priori,
the support of the solution $v(t,.)$ is greater or equal to the support of $%
u_{\lambda ^{\ast }}$ for any $t>0$. The following result gives some
indication about how the support of $v(t,.)$ should increase slowly with
time. We shall apply the general local energy methods for the study of free
boundary problems (see, e.g. \cite{ADS}). Notice that for our goal we only
need to get some information on $v(t,.)$ on the level sets where this
function is small enough. So, given $\theta >0$ and $t\geq 0$ we introduce
the notation%
\begin{equation*}
\Omega _{v,\theta }(t):=\{x\in \Omega :v(t,x)\leq \theta \}.
\end{equation*}

\begin{thm}
\label{free boundary1} Assume $\lambda >\lambda ^{\ast }$, $v_{0}=u_{\lambda
^{\ast }}$ and let $\theta >0$ such that $\theta ^{\beta -\alpha }<1/\lambda
$. Let $x_{0}\in \mathbb{R}^{N}\setminus $
supp($v_{0})$ such that $B_{\rho _{0}}(x_{0})\subset
\mathbb{R}^{N}\setminus $supp($v_{0})$ for some $\rho _{0}>0$. Then there
exists $\widehat{t}>0$ and a continuous decreasing function $\rho :[0,%
\widehat{t}]\rightarrow \lbrack 0,\rho _{0}]$ such that $\rho (0)=\rho _{0}$%
, $\rho (\widehat{t})=0$ and $B_{\rho (t)}(x_{0})\subset \mathbb{R}%
^{N}\setminus $supp($v(t,.)\cap \Omega _{v,\theta }(t))$ for any $t\in
\lbrack 0,\widehat{t}]$. In particular, $v(t,x)=0$ a.e. $x\in B_{\rho
(t)}(x_{0})$ for any $t\in \lbrack 0,\widehat{t}].$
\end{thm}

\noindent {\em Proof}\quad It is enough to apply Theorem 2.2 of \cite{ADS}
to the special case of $\psi (u)=u$  and
\begin{equation*}
A(x,t,u,Du)=Du,~~ B(x,t,u,Du)=0,~~C(x,t,u,Du)=(1-\lambda \theta ^{\beta -\alpha })|u|^{\alpha -1}u,
\end{equation*}%
since we know that
\begin{equation*}
v_{t}-\Delta v+(1-\lambda \theta ^{\beta -\alpha })|v|^{\alpha -1}v\leq 0%
\text{ on }\cup _{t>0}\{t\}\times \Omega _{v,\theta }(t),
\end{equation*}%
and all the assumptions of Theorem 2.2 of \cite{ADS} hold.$\fin$

\bigskip

When $\lambda =\lambda ^{\ast }$ we can also give an idea how the support of
$v(t,.)$ corresponding to a strictly positive initial decreases, after a
finite time large enough (remember that in that case the support of $v(t,.)$
must decrease from $\overline{\Omega }$ to the closed ball of $\mathbb{R}^{N}
$ of radio $R(\Omega )\subset \overline{\Omega }$). In this case, we shall
pay attention to the special choice of $v_{0}=u_{\lambda }$ for some $%
\lambda >\lambda ^{\ast }$. Notice that now $u_{\lambda }$ is a
supersolution to $PP(\alpha ,\beta ,\lambda ^{\ast },v_{0})$ since

\begin{equation*}
(u_{\lambda })_{t}-\Delta u_{\lambda }+|u_{\lambda }|^{\alpha -1}u_{\lambda
}=\lambda |u_{\lambda }|^{\beta -1}u_{\lambda }>\lambda ^{\ast }|u_{\lambda
^{\ast }}|^{\beta -1}u_{\lambda ^{\ast }}.
\end{equation*}%
As above, if $u_{\lambda }$ is nondegenerate, we can even prove that $%
v_{t}\leq 0$ a.e. $(0,+\infty )\times \Omega $. \ Concerning the formation
of the free boundary we have:

\begin{thm}
\label{free boundary2} Assume $\lambda =\lambda ^{\ast }$, $v_{0}=u_{\lambda
}$ for some $\lambda >\lambda ^{\ast }$and let $\theta >0$ such that $\theta
^{\beta -\alpha }<1/\lambda ^{\ast }$. Then, for any time $T>0$ large
enough, there exist a finite time $t^{\#}>0$ and a continuous increasing
function $\rho :[t^{\#},T]\rightarrow \lbrack 0,+\infty )$ such that $\rho
(t^{\#})=0$, and $B_{\rho _{0}}(x_{0})\subset \mathbb{R}^{N}\setminus $%
support($v(t,.)\cap \Omega _{v,\theta }(t))$ for any $t\in \lbrack t^{\#},T]$%
. In particular, $v(t,x)=0$ a.e. $x\in B_{\rho (t)}(x_{0})$ for any $t\in
\lbrack t^{\#},T].$
\end{thm}

\noindent {\em Proof}\quad This time it is enough to apply Theorem 4.2 of
\cite{ADS} to the special case of $\psi (u)=u$, $A(x,t,u,Du)=Du$, $%
B(x,t,u,Du)=0$ and
\begin{equation*}
C(x,t,u,Du)=(1-\lambda ^{\ast }\theta ^{\beta -\alpha })|u|^{\alpha -1}u.
\end{equation*}%
Indeed, as above we know that
\begin{equation*}
v_{t}-\Delta v+(1-\lambda ^{\ast }\theta ^{\beta -\alpha })|v|^{\alpha
-1}v\leq 0\text{ on }\cup _{t\in (0,T)}\{t\}\times \Omega _{v,\theta }(t).
\end{equation*}%
and all the assumptions of Theorem 4.2 of \cite{ADS} hold.$\fin$

\bigskip

\section{Appendix}

\begin{lem}
\label{app} Assume $\lambda \in [\Lambda_{1P}, +\infty)$ and $u_{\lambda_m}$
is a sequence of solutions of (\ref{min1}), where $\lambda_m \to \lambda$ as
$m\to +\infty$. Then there exist a minimizer $u_\lambda$ of (\ref{min1}) and
a subsequence, still denoted by $(u_{\lambda_m})$, such that $u_{\lambda_m}
\to u_\lambda$ strongly in $H_{0}^{1}$ as $m \to +\infty$.
\end{lem}

\noindent {\em Proof}\quad Let $\lambda \in \lbrack \Lambda _{1P},+\infty )$%
, $\lambda _{m}\rightarrow \lambda $ as $m\rightarrow +\infty $ and $%
u_{\lambda _{m}}$ be a sequence of solutions of (\ref{min1}). As in the
proof of Lemma \ref{le1e} it is derived that the set $(u_{\lambda _{m}})$ is
bounded in $H_{0}^{1}$. Hence by the Sobolev embedding theorem there exists
a subsequence, still denoted by $(u_{\lambda _{m}})$, such that
\begin{equation}
u_{\lambda _{m}}\rightharpoondown \bar{u}_{\lambda }~~\mbox{weakly in}%
~~H_{0}^{1},~~~u_{\lambda _{m}}\rightarrow \bar{u}_{\lambda }~~%
\mbox{strongly in}~~L_{q}(\Omega ),  \label{convAp}
\end{equation}%
where $0<q<2^{\ast }$, for some limit point $\bar{u}_{\lambda }$. As in the
proof of Lemma \ref{le1e} one derives that $\bar{u}_{\lambda }\neq 0$ and
\begin{equation}\label{C1}
E_{\lambda }(\bar{u}_{\lambda })\leq \liminf_{m\rightarrow \infty
}E_{\lambda _{m}}(u_{\lambda _{m}}),~~E_{\lambda }^{\prime }(\bar{u}%
_{\lambda })\leq 0,~~P_{\lambda }(\bar{u}_{\lambda })\leq 0.
\end{equation}%
Let $\lambda >\Lambda _{1P}$. By Lemma \ref{le1e} there exists a minimizer $%
u_{\lambda }$ of (\ref{min1}), i.e. $u_{\lambda }\in M_{\lambda }$ and $\hat{%
E}_{\lambda }=E_{\lambda }(u_{\lambda })$. Then
\begin{equation*}
|E_{\lambda }(u_{\lambda })-E_{\lambda _{m}}(u_{\lambda })|<C|\lambda
-\lambda _{m}|,
\end{equation*}%
where $C<+\infty $ does not depend on $m$. Furthermore,
\begin{equation*}
E_{\lambda _{m}}(u_{\lambda })\geq E_{\lambda _{m}}(t_{\min }(u_{\lambda
})u_{\lambda })\geq E_{\lambda _{m}}(u_{\lambda _{m}})
\end{equation*}%
provided that $m$ is a sufficiently large number. Thus  we have
\begin{equation*}
E_{\lambda }(u_{\lambda })+C|\lambda -\lambda _{m}|>E_{\lambda
_{m}}(u_{\lambda })\geq E_{\lambda _{m}}(u_{\lambda _{m}}),
\end{equation*}%
and therefore $\hat{E}_{\lambda }:=E_{\lambda }(u_{\lambda })\geq
\liminf_{m\rightarrow \infty }E_{\lambda _{m}}(u_{\lambda _{m}})$. Hence by %
\eqref{C1} we have
\begin{equation*}
E_{\lambda }(\bar{u}_{\lambda })\leq \hat{E}_{\lambda }.
\end{equation*}%
Assume $E_{\lambda }^{\prime }(\bar{u}_{\lambda })<0$. Then $E_{\lambda
}^{\prime }(t_{\min }(\bar{u}_{\lambda })\bar{u}_{\lambda })=0$ and $%
E_{\lambda }(t_{\min }(\bar{u}_{\lambda })\bar{u}_{\lambda })<E_{\lambda }(%
\bar{u}_{\lambda })\leq \hat{E}_{\lambda }$. In virtue of Proposition \ref%
{pradd}, this implies that $P_{\lambda }(t_{\min }(\bar{u}_{\lambda })\bar{u}%
_{\lambda })<0$. Thus $t_{\min }(\bar{u}_{\lambda })\bar{u}_{\lambda }\in
M_{\lambda }$ and since $E_{\lambda }(t_{\min }(\bar{u}_{\lambda })\bar{u}%
_{\lambda })<\hat{E}_{\lambda }$ we get a contradiction. Hence $E_{\lambda }(%
\bar{u}_{\lambda })=\hat{E}_{\lambda }$, $E_{\lambda }^{\prime }(\bar{u}%
_{\lambda })=0$ and $u_{\lambda _{m}}\rightarrow \bar{u}_{\lambda }$
strongly in $H_{0}^{1}$ as $m\rightarrow +\infty $ .

Assume now that $\lambda =\Lambda _{1P}$. Since $E_{\lambda _{m}}^{\prime
}(u_{\lambda _{m}})=0$, $P_{\lambda _{m}}(u_{\lambda _{m}})\leq 0$, we have $%
r_{u_{\lambda _{m}}}^{P}(1)\leq \lambda _{m}=r_{u_{\lambda
_{m}}}^{1}(1)$. Then by Proposition \ref{propEst} (see Figure 5),
$1\in \lbrack t_{1P}(u_{\lambda _{m}}),+\infty )$ and therefore
\begin{equation*}
\lambda _{1P}(u_{\lambda _{m}})=r_{u_{\lambda _{m}}}^{1}(t_{1P}(u_{\lambda
_{m}}))\leq r_{u_{\lambda _{m}}}^{1}(1)=\lambda _{m},~~m=1,2,....
\end{equation*}%
Hence, since $\lambda _{m}\downarrow \lambda $, we have $\lambda
_{1P}(u_{\lambda _{m}})\downarrow \Lambda _{1P}$ as $m\rightarrow \infty $.
Thus, $(u_{\lambda _{m}})$ is a minimizing sequence of \eqref{PPoh} and
therefore by \eqref{convAp}, $\lambda _{1P}(\bar{u}_{\lambda })\leq \Lambda
_{1P}$. Since the strong inequality $\lambda _{1P}(\bar{u}_{\lambda
})<\Lambda _{1P}$ is impossible, we conclude that $\lambda _{1P}(\bar{u}%
_{\lambda })=\Lambda _{1P}$, which yields that $u_{\lambda _{m}}\rightarrow
\bar{u}_{\lambda }$ strongly in $H_{0}^{1}.\fin$ 
\bigskip

\textbf{Acknowledgments}

The research of J.I. D\'{\i}az and J. Hern\'{a}ndez was partially supported
by the projects ref. MTM 2014-57113-P and MTM2017-85449-P of the DGISPI
(Spain).

\flushright{\small
\begin{tabular}{llll}
J.I. D\'{\i}az & J.~Hern\'{a}ndez & Y.Sh.~Ilyasov \\
Instituto de Matem\'{a}tica Interdisciplinar & ,Instituto de Matem\'{a}tica Interdisciplinar & 
Institute of Mathematics of UFRC RAS \\
Universidad Complutense de Madrid& Universidad Complutense de Madrid & Chernyshevsky str., 450008, Ufa, Rusia\\
28040 Madrid, Spain & 28040 Madrid, Spain & Instituto de Matem\'atica e Estat\'istica.\\
& &  Universidade Federal de Goi\'as,\\
& & 74001-970, Goiania, Brazil\\
{\tt jidiaz@ucm.es} & {\tt jesus.hernande@telefonica.net}  & {\tt ilyasov02@gmail.com}
\end{tabular}
}


\begin{thebibliography}{99}
\bibitem{Annello} A. Anello and F. Faraci, Two solutions for an elliptic
problem with two singular terms. \textit{Calc. Var. and PDEs}
\textbf{56} (2017), 31 pages.

\bibitem{ADS} S. Antontsev, J. I. D\'{\i}az, S. Shmarev, \textit{Energy
methods for free boundary problems. Applications to nonlinear PDEs
and Fluid Mechanics}, Birk\"{a}user, Boston, 2002.

\bibitem{CDE} T. Cazenave, T. Dickstein and M. Escobedo, A semilinear heat
equation with concave-convex nonlinearity, \textit{Rendiconti di Matematica}%
, \textit{Serie VII}, \textbf{19} (1999) 211-242.

\bibitem{CortElgFelmer-1} C.~Cort\'{a}zar, M.~Elgueta, and P.~Felmer,
Symmetry in an elliptic problem and the blow-up set of a quasilinear
heat equation, \textit{Comm. P.D.E.}, \textbf{21} (1996) 507-520.

\bibitem{CortElgFelmer-2} C.~Cort\'{a}zar, M.~Elgueta, and P.~Felmer, On a
semi-linear elliptic problem in $\mathbb{R}^{N}$ with a
non-Lipschitzian non\--linearity, \textit{Advances in Diff. Eqs.},
\textbf{1} (1996) 199-218.

\bibitem{diaz} J.I. D\'{\i}az and J. Hern\'{a}ndez,, Global bifurcation and
continua of nonnegative solutions for a quasilinear elliptic
problem, \textit{C.R. Acad. Sci. Paris}, \textbf{329}, (1999),
587-592.

\bibitem{Diaz-Hernan-Man} J.~I. D{\'{\i}}az, J. Hern\'{a}ndez and F.~J.
Mancebo, Branches of positive and free boundary solutions for some
singular quasi\-linear elliptic problems, \textit{J. Math. Anal.
Appl.}, \textbf{352} (2009), 449--474.

\bibitem{DIH1} J. I. D\'{\i}az, J. Hern\'{a}ndez and \ Y\textit{. }Sh.
Il'yasov, On the existence of positive solutions and solutions with
compact support for a spectral nonlinear elliptic problem with
strong absorption, \textit{Nonl. Anal.: Th., Meth. \& Appl.}
\textbf{119} (2015), 484--500.

\bibitem{DIH} J. I. D\'{\i}az, J. Hern\'{a}ndez and Y\textit{. }Sh.
Il'yasov, Stability criteria on flat and compactly supported ground
states of some non-Lipschitz autonomous semilinear equations,
\textit{Chinese Ann. Math.} \textbf{38} (2017), 345-378.

\bibitem{GazzolaSerin} F. Gazzola, J. Serrin and M. Tang, Existence of
ground states and free boundary problems for quasilinear elliptic
operators, \textit{Advances in Differential Equations} \textbf{5}
(2000) 1-30.

\bibitem{giltrud} D. Gilbarg and N.S. Trudinger, \textit{Elliptic partial
differential equations of second oder}. 2nd edition, Grundlehren
224, Springer, Berlin-Heidelberg-New York-Tokyo 1983.

%\bibitem{ilyas} Il'yasov, Y.S., \textit{Nonlocal investigations of
%bifurcations of solutions of nonlinear elliptic equations}, Izv. Math. 66
%no. 6,(2002), 1103--1130.

\bibitem{HMV} J. Hern\'{a}ndez, F. Mancebo and J. M. Vega, Positive
solutions for singular nonlinear elliptic equations, \textit{Proc.
Roy. Soc. Edinburgh} \textbf{137A} (2007), 41-62.

\bibitem{ilDr} Y. Sh. Il'yasov, On critical exponent for an elliptic
equation with non-Lipschitz nonlinearity, \textit{Dynamical
Systems,} Supplement (2011), 698-706.

\bibitem{ilCrit} Y. Sh. Il'yasov, On the curve of critical exponents for
nonlinear elliptic problems in the case of a zero mass. \textit{%
Computational Mathematics and Mathematical Physics} \textbf{57} 3
(2017), 497-514.

\bibitem{IlEg} Y. Sh. Ilyasov and Y. Egorov, Hopf maximum principle
violation for elliptic equations with non-Lipschitz nonlinearity, \textit{%
Nonlin. Anal.} \textbf{72} (2010) 3346-3355.

%\bibitem{ilegorov} Y. Il'yasov, Y. Egorov, Hopf boundary maximum principle violation for semilinear elliptic equations, Nonlinear Analysis: Theory, Methods \& Applications 72 (7-8), 3346-3355, (2010)

\bibitem{ilyaReil} Y. Sh. Ilyasov, On extreme values of Nehari manifold
method via nonlinear Rayleigh's quotient. \textit{Topological
Methods in Nonlinear Analysis} \textbf{49} 2 (2017) 683-714.

\bibitem{Takac_Ilyasov} Y. Sh. Il'yasov and P. Takac, Optimal-regularity,
Pohozhaev's identity, and nonexistence of weak solutions to some
quasilinear elliptic equations. \textit{Journal of Differential
Equations}, \textbf{252} 3 (2012) 2792-2822.

\bibitem{Kaper1} H. Kaper and M. Kwong, Free boundary problems for
Emden-Fowler equation, \textit{Differential and Integral Equations}, \textbf{%
3} (1990) 353-362.

\bibitem{Kaper2} H. Kaper, M. Kwong and Y. Li, Symmetry results for
reaction-diffusion equations, \textit{Differential and Integral
Equations}, \textbf{6} (1993) 1045-1056.

\bibitem{Lieberman} G.~M. Lieberman, Boundary regularity for solutions of
degenerate elliptic equations, \textit{Nonlinear Anal.}, \textbf{12
}11 (1988) 1203--1219.

\bibitem{Montenegro} M. Montenegro and E. Silva, Two solutions for a
singular elliptic equation by variational methods. \textit{Ann. Sc.
Norm. Sup. Pisa} \textbf{11 }(2012), 143-165.

\bibitem{poh} S.I. Pohozaev, Eigenfunctions of the equation $\Delta
u+\lambda f(u)=0$, \textit{Sov. Math. Doklady} \textbf{5} (1965)
1408-1411.

\bibitem{PucciSerrin} P. Pucci and J. Serrin, Uniqueness of ground states
for quasilinear elliptic operators. \textit{Indiana University
Mathematics
Journal} \textbf{47} 2 (1998) 501-528. %\bibitem{reich} Reichel, W.,
%Uniqueness theorems for variational problems by the method of transformation groups.
%Lecture Notes in Mathematics 1841. Berlin: Springer (2004).

\bibitem{Rosenau} P. Rosenau and J. M. Hyman, Compactons: Solitons with
finite wavelength, \textit{Phys. Rev. Lett.} \textbf{70} 5 (1993)
564--567.

\bibitem{Serrin-Zou} J. Serrin and H. \ Zou, Symmetry of ground states of
quasilinear elliptic equations. \textit{Archive for Rational
Mechanics and Analysis}, \textbf{148} 4 (1999) 265-290.

\bibitem{str} M. Struwe, \textit{Variational Methods, Application to
Nonlinear Partial Differential Equations and Hamiltonian Systems}.
Springer- Verlag, Berlin, Heidelberg, New-York, 1996.

\bibitem{Temam} R.~Temam, \textit{Navier-Stokes Equations: Theory and
Numerical Analysis}, North-Holland, Amsterdam--New York--Oxford,
1979.

\bibitem{Zeidler} E. Zeidler, \textit{Nonlinear functional analysis,} Vol.3.
\textit{Variational methods and optimization}, Springer-Verlag,
Berlin, 1985.

\end{thebibliography}
\end{document}